\documentclass[a4paper,11pt]{article}
\usepackage[latin1]{inputenc}
\usepackage{amsfonts}
\usepackage{amsmath}
\usepackage{amssymb}
\usepackage{amsthm}
\usepackage[T1]{fontenc}
\usepackage[dvips]{graphicx}
\usepackage{color}
\usepackage{epsfig}
\theoremstyle{plain}
\newtheorem{Thm}{Theorem}
\newtheorem{Lem}{Lemma}
\newtheorem{Def}{Definition}
\newtheorem{Cor}{Corollary}
\newtheorem*{Rem*}{\textsc{Remark}}
\newtheorem*{Lem*}{\textsc{Lemma}}
\newtheorem*{Cor*}{\textsc{Corollary}}
\newtheorem*{Con*}{\textsc{Conjecture}}
\newcommand{\e}{\varepsilon}
\newcommand{\av}[1]{\left| \! #1 \! \right|}

\title{Poles of Int\'egrale Tritronqu\'ee and Anharmonic Oscillators. A WKB Approach}
\author{Davide Masoero \thanks{E-mail address: masoero@sissa.it}\\ SISSA - Trieste}
\date{}
\begin{document}

\maketitle

\begin{abstract}
Poles of solutions to the Painlev\'e-I equations are intimately related to the
theory of the cubic anharmonic oscillator.
In particular, poles of int\'egrale tritronqu\'ee are in bijection with cubic oscillators
that admit the simultaneous solutions of two quantization conditions.
We analyze this pair of quantization conditions by means of a suitable version of the
complex WKB method.
\end{abstract}

\section{Introduction}
The aim of the present paper is to study the distribution of poles
of the solutions $y=y(z)$
to the Painlev\'e first equation (P-I) 

\begin{equation*}
y''= 6y^2 -z \, , \; z \in \mathbb{C} \quad ,
\end{equation*}
with particular attention to the poles of the \textit{int\'egrale tritronqu\'ee}.

As it is well-known, any local solution of P-I
extends to a global meromorphic function $y(z), z\in \mathbb{C}$, with
an essential singularity at infinity \cite{gromak}. Global solutions of P-I are
called Painlev\'e-I transcendents, since they cannot
be expressed via elementary functions or classical special functions \cite{ince}.
The int\'egrale tritronqu\'ee is a special P-I transcendent, which was
discovered by
Boutroux in his classical paper \cite{boutroux}
(see \cite{joshikruskal} and \cite{elliptic} for a modern review).
Boutroux characterized the int\'egrale tritronqu\'ee as the unique solution
of P-I with the following asymptotic behaviour at infinity

\begin{equation*}
y(z) \sim - \sqrt{\frac{z}{6}}, \quad \mbox{if} \quad |\arg z| <\frac{4 \pi}{5} \; .
\end{equation*}

Nowadays Painlev\'e first equation is studied in many areas of mathematics and physics.
Indeed, it is remarkable that special solutions of P-I describe semiclassical
asymptotics of a wealth of different problems (see \cite{kapaev} and references therein).
In particular, in \cite{dubrovin}
it is discovered that the int\'egrale tritronqu\'ee provides
the universal correction to the dispersionless limit of solutions to the focusing nonlinear
Schr\"odinger equation.

Theoretical and numerical evidences led the authors of \cite{dubrovin}
to the following inspiring 

\begin{Con*}
If $a \in \mathbb{C}$ is a pole of the int\'egrale tritronqu\'ee then $\av{\,\arg{a}\,}\! \geq\! \frac{4 \pi}{5}$.
\end{Con*}

Following the \textit{isomonodromic} approach to P-I \cite{kapaev}, any
solution  $y(z)$  gives rise to an
\textit{isomonodromic deformation}
of the following linear equation with an irregular singularity

\begin{eqnarray*}
\overrightarrow{\Phi}_{\lambda}(\lambda,z)&=&
\left(
\begin{matrix}
y'(z) & 2\lambda^2+2\lambda y(z)-z+2y^2(z) \\ 
2(\lambda-y(z)) & -y'(z)
\end{matrix} \right)
 \overrightarrow{\Phi}(\lambda,z)  \; .
\end{eqnarray*}

The deformation of the equation is manifestly singular at every pole $a \in \mathbb{C}$ of $y$,
however in Theorem \ref{them:toCubic} we show
that at the singularity this equation  can be replaced with a simpler one,
which has the same \textit{monodromy data} (cf. \cite{its} for Painlev\'e II).
This is the following  Schr\"odinger equation with cubic potential

\begin{equation*}
\frac{d^2\psi(\lambda)}{d\lambda^2}= V(\lambda;a,b) \psi(\lambda)\; , \quad
V(\lambda;a,b)=4 \lambda^3 -2 a \lambda -28 b \;.
\end{equation*}
Here $a$ is the location of the pole of $y$ and $b$ is a complex number entering into the
Laurent expansion of $y$ around $a$ (see formula (\ref{laurent}) below).

The isomonodromy
property implies that there exists a natural
injective map $\cal{M}$ from the space of solutions of P-I
to the space of \textit{monodromy data} of the above equations
(see Lemma \ref{lem:multipliers}), while
the Schr\"odinger equation defines naturally a map
$\cal{T}$ from the space of cubic potentials
to the space of monodromy data.

Our first main result is Theorem \ref{lem:bijection}, which
states that  $a \in \mathbb{C}$ is a pole of $y(z)$ \textit{if and only if}
there exists $b \in \mathbb{C}$ such that
${\cal{M}}(y)={\cal{T}}(V(\lambda;a,b))$.

In particular, because of the special monodromy data related to the int\'egrale tritronqu\'ee
(see  Theorem \ref{the:kap}, due to Kapaev), we will show that
the poles of the int\'egrale tritronqu\'ee  are
in bijection with the
simultaneous solutions of two different
\textit{quantization conditions}.

The above approach naturally embeds the study of poles
of Painlev\'e-I transcendents into the Nevanlinna's theory of branched coverings
of the sphere and the complex WKB method.

The beautiful theory of R. Nevanlinna (see  \cite{nevanlinna} and \cite{elfving})
relates bijectively the Schr\"odinger equations with a polynomial potential to the
branched coverings of the sphere with logarithmic branch points, considered up to conformal equivalence. 
Using this theory we are able to prove
the surjectivity of the map $\cal{M}$ (see Theorem \ref{thm:surjection}).

Moreover Nevanlinna's theory provides the poles of any solution of P-I
with an unexpected and remarkable rich structure.
In particular, poles of the tritronqu\'ee solution
can be labelled by the monodromy of coverings of the Riemann sphere
with 3 logarithmic branch points. In a subsequent paper,
we are going to use this topological description to complete the WKB analysis
of the present paper.

The WKB analysis of P-I developed in \cite{takei}
has never been applied to the direct study of the distributions of poles.
To achieve such a goal we follow the Fedoryuk's approach (see \cite{fedoryuk})
to the complex WKB theory,
and in the Classification Theorem  we give a complete topological classification
of the \textit{Stokes complexes} of all  cubic potentials.
As a consequence of the Classification Theorem,
we obtain our second main result:
all polynomials whose monodromy data, in the WKB approximation,
are the monodromy data of the int\'egrale tritronqu\'ee 
have the same topological type of Stokes complex
and satisfy a pair of \textit{Bohr-Sommerfeld} quantization conditions, namely system (\ref{eq:boutroux}).
In particular, in this way we reproduce the conditions obtained by Boutroux,
through a completely different approach,
in his study of the asymptotic distributions of the poles of the int\'egrale tritronqu\'ee.

A priori, the WKB method is expected to give an approximation
of poles $z=a$ for $a$ sufficiently large. Surprisingly
our approach proves to be numerically very efficient
also for poles close to the origin, see Table 2 below.
\footnote{In the subsequent paper \cite{piwkb2} the author shows that eventually
around any solution of the Bohr-Sommerfeld-Boutroux system there is one and only one pole
of the int\'egrale tritronqu\'ee and the distance between a pole and its approximation
vanishes asymptotically.}

The paper is organized as follows.
In Section 2 we derive the Schr\"odinger equation associate with P-I and study
thoroughly its relations with poles
of P-I transcendents. Section 3 is devoted to the topological classification of Stokes complexes. 
In Section 4 we calculate the monodromy data in the WKB approximation, we derive
the correct Bohr-Sommerfeld conditions for the poles of tritronqu\'ee, and
we introduce the "small parameter" of the approximation.
In Section 5 we obtain an asymptotic description of poles of the \textit{integrale tritronqu\'ee}.
In Appendix A and Appendix B we prove some theorems regarding the WKB
functions that are used in section 2 and 3.

\paragraph{Acknowledgments}
I am indebted to Prof. B. Dubrovin who introduced me to the problem and constantly gave me
suggestions and advice.
I would like to thank A. Raimondo, G. De Nittis and R. Tateo for useful discussions.

Part of this work was written during a visit
at the department "R. Caccioppoli" of Napoli University.
This work is partially supported by the European Science Foundation Programme
"Methods of Integrable Systems, Geometry, Applied Mathematics" (MISGAM), and
by the Italian Ministry of University and Research (MIUR) grant PRIN 2006
"Geometric methods in the theory of nonlinear waves and their applications".

\section{Poles and Cubic Oscillators}

We review some well-known facts regarding the isomonodromic approach to the P-I equation
and analyze the isomonodromic deformation in a neighborhood of the singularities.

\subsection{P-I as an Isomonodromic Deformation}
P-I is equivalent to the compatibility condition of the following system of linear ODEs:

\begin{eqnarray}\label{sys:1}
\overrightarrow{\Phi}_{\lambda}(\lambda,z)&=&
\left(
\begin{matrix}
y'(z) & 2\lambda^2+2\lambda y(z)-z+2y^2(z) \\ 
2(\lambda-y(z)) & -y'(z)
\end{matrix} \right)
 \overrightarrow{\Phi}(\lambda,z) \\ \label{sys:2}
\overrightarrow{\Phi}_{z}(\lambda,z)&=&-
\left(
\begin{matrix}
0 & 2y(z)+\lambda \\ 
1  & 0
\end{matrix} \right)
 \overrightarrow{\Phi}(\lambda,z) \quad .
\end{eqnarray}

The precise meaning of the word \textit{compatibility} is given by the following

\begin{Lem}\label{lem:compatibility}
Fix $z_0$, $\lambda_0$ and the Cauchy data $y(z_0)$, $y'(z_0)$, and $\overrightarrow{\Phi}(\lambda_0, z_0)$.
Let $U_{z_0}$ be any simply connected neighborhood of $z_0$.
Then $y(z)$ satisfies the  Painlev\'e first equation in  $U_{z_0}$ \textit{iff}
the system (\ref{sys:1},\ref{sys:2}) has a solution $\forall (\lambda,z) \in \mathbb{C} \times U_{z_0}$.
Moreover the solution is unique.
\begin{proof}
See \cite{kapaev}.
\end{proof}
\end{Lem}

In this subsection we suppose that we have fixed a solution $y$ of P-I and
a simply connected region $U$ such that $y|_U$ is holomorphic.

We are now going to define the important concepts of \textit{monodromy data}
and \textit{isomonodromic deformation}
of equation (\ref{sys:1}). For this reason, we have to introduce some
particular solutions
of system (\ref{sys:1},\ref{sys:2}), to be
uniquely defined by the asymptotic behaviour for $\lambda \to \infty$.

Fix $k \in \mathbb{Z}_5 = \left\lbrace -2, \dots,2 \right\rbrace$
and the branch of $\lambda^{\frac{1}{2}}$ in such a way that
${\rm Re}{\lambda^{\frac{5}{2}}} \to + \infty$ as
$|\lambda| \to \infty,\arg{\lambda}=\frac{2\pi k}{5}$.
Then (see \cite{kapaev}) for any $y$ solution of P-I, there exists a unique  solution
$\overrightarrow{\Phi_k}(\lambda,z)$
of (\ref{sys:1},\ref{sys:2})
such that 
\begin{equation} \label{asymptotic}
 \lim_{\substack{\lambda  \to \infty \\ \left| \arg{\lambda}-2\frac{\pi k}{5} \right| <
\frac{3 \pi}{5}-\varepsilon}}
e^{+\frac{4}{5}\lambda^{\frac{5}{2}} -\frac{1}{2} z\lambda^{\frac{1}{2}}} \left(
\begin{matrix}
\lambda^{-\frac{1}{4}} & 0 \\ 
0  & \lambda^{+\frac{1}{4}}
\end{matrix} \right) \overrightarrow{\Phi}_k(\lambda)= \left(
\begin{matrix}
1 \\ 
1
\end{matrix} \right), \forall z \in U \; ,
\end{equation}
where  $\lambda^{\frac{1}{4}}$ is defined globally on the complex plane minus the negative
real axis, and is positive on the positive real axis. Notice that, depending on $k \in \mathbb{Z}_5$, $\left(\lambda^{\frac{1}{4}}\right)^2$ may not be equal to $\lambda^{\frac{1}{2}}$.
Here and in the following, if not otherwise stated, $\varepsilon$ is an arbitrarily small positive number.

From the asymptotics (\ref{asymptotic}) it follows that $\overrightarrow{\Phi}_k(\lambda,z)$ and
$\overrightarrow{\Phi}_{k+1}(\lambda,z)$ are linearly independent for
any $k \in \mathbb{Z}_5$ and the following equality
holds true
\begin{equation}\label{eq:sigmak}
\overrightarrow{\Phi}_{k-1}(\lambda)=\overrightarrow{\Phi}_{k+1}(\lambda)+\sigma_k(z)
\overrightarrow{\Phi}_{k}(\lambda) \; ,
\end{equation}
where $\sigma_k(z)$ is an analytic function of $z$, for any $k \in \mathbb{Z}_5$.

\begin{Def}
Fixed $z$, we call $\sigma_k$ the k-th {\rm Stokes multiplier} of equation (\ref{sys:1})
and the set of all five Stokes multipliers the {\rm monodromy data} of
(\ref{sys:1}). The problem of calculating the monodromy data is called the {\rm direct monodromy problem}.
\end{Def}

Stokes multipliers are very important for our analysis and we list their main properties in the following

\begin{Lem}\label{lem:multipliers}
Let $\sigma_k(z), k\in \mathbb{Z}_5$ be defined as above. Then
\begin{itemize}
 \item[(i)] equation (\ref{sys:2}) is an \textit{isomonodromic deformation} of equation (\ref{sys:1}),
 i.e. $\frac{d \sigma_k(z)}{dz}=0$.
 \item[(ii)] The numbers $ \sigma_k, k \in \mathbb{Z}_5$ satisfy the following system of algebraic equations 
\begin{equation} \label{eq:quadratic}
1+ \sigma_k\sigma_{k+1}= -i \, \sigma_{k+3} \, , \; k \in \mathbb{Z}_5 \;.
       \end{equation}
\end{itemize}
\begin{proof}
\begin{itemize}
See \cite{kapaev}.
\end{itemize}

\end{proof}

\end{Lem}
Observe that only 3 of the algebraic equations (\ref{eq:quadratic}) are independent.
\begin{Def}
We denote $V$ the algebraic variety of quintuplets of complex numbers
satisfying (\ref{eq:quadratic}) and call \textit{admissible monodromy data}
the elements of $V$.
Due to Lemma \ref{lem:multipliers}, equations (\ref{eq:sigmak}) define the following map
\begin{equation*}
{\cal{M}}: \left\lbrace \mbox{P-I transcendents} \right\rbrace \to V \; .
\end{equation*}
\end{Def}

\begin{Lem*}
 $\cal{M}$ is injective.
\begin{proof}
  See \cite{kapaev}.
\end{proof}
\end{Lem*}

We end the section with a result of Kapaev, which completely characterizes
the int\'egrale tritronqu\'ee in term of Stokes multipliers.

\begin{Thm}\label{the:kap}(Kapaev)

The image under $\cal{M}$ of the int\'egrale tritronqu\'ee are the monodromy data
uniquely characterized by the following equalities

\begin{equation}
\sigma_2=\sigma_{-2}=0 \; .
\end{equation}
\begin{proof}
See \cite{kapaev}.
\end{proof}

\end{Thm}

\subsection{Poles of $y$: cubic oscillator}

So far we dealt with the system (\ref{sys:1},\ref{sys:2}) in a region $U$ which does not contain
any pole of $y(z)$. Indeed, the situation at a pole is different, for equation
(\ref{sys:1}) makes no sense.

However, we show that any solution $\overrightarrow{\Phi}(\lambda,z)$
of system (\ref{sys:1}, \ref{sys:2}) is meromorphic in all the $z-plane$;  moreover,
a pole $z=a$ of $y(z)$ is also a pole of $\overrightarrow{\Phi}(\lambda,z)$ and the residue at the pole
of its second component satisfies the scalar equation of Schr\"odinger type (\ref{eq:schr}).

In order to be able to describe the local behavior of
$\overrightarrow{\Phi}(\lambda,z)$ near a pole $a$ of $y(z)$, we have to know
the local behavior of $y(z)$ close to the same point $a$.

\begin{Lem}[Painlev\'e]\label{lem:laurent}
Let $a \in \mathbb{C}$ be a pole of $y$. Then in a neighborhood of ${\rm a}$, $y$
has the following convergent Laurent expansion
\begin{equation}\label{laurent}
y(z)=\frac{1}{(z-a)^2}+\frac{a(z-a)^2}{10}+ \frac{(z-a)^3}{6}+b(z-a)^4+\sum_{j \geq 5}c_j(a,b) (z-a)^j
\end{equation}
where ${b}$ is some complex number and $c_j(a,b)$ are real polynomials in ${a}$ and ${b}$, not depending
on the particular solution $y$.

Conversely, fixed arbitrary $a,b \in \mathbb{C}$, the above expansion has a non zero
radius of convergence and solves P-I.

\begin{proof}
 See \cite{gromak}.
\end{proof}
\end{Lem}

\begin{Def}
We define the  map 
\begin{equation*}
{\cal{L}}: \mathbb{C}^2 \to \left\lbrace \mbox{P-I transcendents} \right\rbrace  \; .
\end{equation*}
${\cal{L}}(a,b)$ is the unique analytic continuation of the  Laurent expansion (\ref{laurent}).
\end{Def}

We have already collected all elements necessary to formulate the important

\begin{Thm}\label{them:toCubic}
Fix a solution $y$ of P-I and let  $\Phi_k^{(i)}(\lambda,z), \; i=1,2 \, \; k \in \mathbb{Z}_5 \;$ be 
the i-th component of $\overrightarrow{\Phi}_k(\lambda,z)$. Then
\begin{itemize}
 \item[(i)]$\overrightarrow{\Phi}_k(\lambda,z)$ is a meromorphic function of $z$.
  All the singularities are double poles. Moreover, $a \in \mathbb{C}$ is a pole of $\overrightarrow{\Phi}_k(\lambda,z)$
  \textit{iff} it is a pole of $y$.
 \item[(ii)] If $a \in \mathbb{C}$ is a pole of $y$ then 
 \begin{equation*}
 \Psi_k(\lambda)=\lim_{z\to a}\left(z-a\right)\Phi_k^{(2)}(\lambda,z)
 \end{equation*}
 is an entire function of $\lambda$. It satisfies
 the following \textit{Schr\"odinger equation with cubic potential}
\begin{equation}\label{eq:schr}
\frac{d^2\Psi_k(\lambda)}{d\lambda^2}= \left( 4 \lambda^3 -2 a \lambda -28 b\right)\Psi_k(\lambda) \; , 
\end{equation}
where $b \in \mathbb{C}$ is the coefficient entering into the Laurent expansion (\ref{laurent}) of $y$ around $a$.
 \item[(iii)]If $\lambda^{\frac{1}{2}}$ and $\lambda^{\frac{1}{4}}$ are chosen as in asymptotics (\ref{asymptotic}), then $\forall \e >0$
\begin{equation}
\lim_{\lambda \to \infty, \av{\lambda -\frac{2 \pi k}{5}} < \frac{3 \pi }{5} -\varepsilon}
\lambda^{\frac{3}{4}}e^{+\frac{4}{5}\lambda^{\frac{5}{2}} - \frac{1}{2}
a\lambda^{\frac{1}{2}}}\Psi_k(\lambda)=i .
\end{equation} 
\item[(iv)]Equation (\ref{eq:schr}) possesses the same monodromy data as equation (\ref{sys:1}), i.e.
\begin{equation*}
 \Psi_{k-1}(\lambda)=\Psi_{k+1}(\lambda)+\sigma_k\Psi_k(\lambda) \; . 
\end{equation*}
\end{itemize}

\begin{proof}

 (i) From the Laurent expansion (\ref{laurent}), it is easily seen that a pole $a$ of $y$ is a fuchsian
 singularity with trivial monodromy of equation (\ref{sys:2}). In particular the following Laurent expansions 
 of $\Phi_k^{(i)}(\lambda,z)$ are valid
\begin{eqnarray}\nonumber
 \! \Phi_k^{(2)}(\lambda,z) \! &=&\! \frac{\psi_k(\lambda)}{(z-a)} \left( 1 \!-\!\frac{\lambda}{2}(z-a)^2\right)
\!+\! \varphi_k(\lambda)(z-a)^2 \!+\! O((z-a)^3) \; ,\\ \label{exp:phi}
\! \Phi_k^{(1)}(\lambda,z) \!&=&\!  \frac{\psi_k(\lambda)}{(z-a)^2} \!
 \left( 1 \!+\!\frac{\lambda}{2} (z-a)^2 \right)\! -\! 2 \varphi_k(\lambda)(z-a) \!+\!  O((z-a)^2) .
\end{eqnarray} 
Expansions (\ref{exp:phi}) show that 
$\overrightarrow{\Phi}_k(\lambda,z)$ is meromorphic in a neighborhood of the point $a$ and  this point is a pole
of order not greater than 2.

(ii), (iii) The proof is in Appendix \ref{app:B}.

(iv)  Since the functions
$(z-a)\Phi_k^{(2)}$ satisfy equations (\ref{eq:sigmak}) for any $z$ with constant Stokes multipliers, then
their limit, i.e. the functions $\Psi_k(\lambda)$,  satisfy the same equations.

\end{proof}
\end{Thm}

\begin{Def}
We call any cubic polynomial of the form 
$V(\lambda; a,b)= 4 \lambda^3 -2 a \lambda -28 b$ a {\rm cubic potential}.
The above formula identifies the space of cubic potentials with $\mathbb{C}^2\backepsilon (a,b)$.

We define the map 
\begin{equation*}
{\cal{T}}: \mathbb{C}^2 \to V  \; .
\end{equation*}
${\cal{T}}(a,b)$ is the monodromy data of equation (\ref{eq:schr}).
\end{Def}

Theorem \ref{them:toCubic} has the following

\begin{Cor*}
$\cal{M}\circ \cal{L}$=$\cal{T}$.
\end{Cor*}

The above corollary implies

\begin{Thm}\label{lem:bijection}
 Let $y$ be any solution of P-I. Then $a \in \mathbb{C}$ is a pole of $y$ iff there exists $b \in \mathbb{C}$
such that ${\cal{M}}(y)={\cal{T}}(a,b)$.
\end{Thm}

We finish this section with a theorem from Nevanlinna's theory \cite{nevanlinna}, which implies the
surjectivity of the map $\cal{M}$.

\begin{Thm}\label{thm:nevanlinna}
The map $\cal{T}$ is surjective. The preimage of any admissible monodromy
data is a countable infinite subset of the space of cubic potentials.
\begin{proof}
See \cite{elfving}.
\end{proof}
\end{Thm}

As a consequence of the above theorem  we have

\begin{Thm}[stated in \cite{kk}] \label{thm:surjection}
The map $\cal{M}$ is bijective: solutions of P-I are in 1-to-1 correspondence with
admissible monodromy data.
\end{Thm}

Theorem \ref{lem:bijection} shows that the distribution of poles of P-I transcendents
is a part of the theory of \textit{anharmonic oscillators},
which has been object of intense study since the seminal papers \cite{wu}
and \cite{simon}.

\begin{Rem*}
In the theory of anharmonic oscillators a special importance is given
to the vanishing of some Stokes multipliers.
For a given $k \in \mathbb{Z}_5$, the
problem is to
find all $(a,b) \in \mathbb{C}^2$ such that the Stokes multiplier
$\sigma_k$ of
equation (\ref{eq:schr}) vanishes. This is called the k-th {\rm lateral connection problem}.
Since fixed $a$, there exists a discrete number
of solutions to any lateral connection problem,
equation $\sigma_k=0$ is referred to as a \textit{quantization condition}.
\end{Rem*}

As a consequence of Theorem \ref{the:kap} and Theorem \ref{lem:bijection}, we have the following
\begin{Cor*}
The point $a \in \mathbb{C}$ is a pole of the int\'egrale tritronqu\'ee if and only if
there exists $b \in \mathbb{C}$ such that the Schr\"odinger equation with the cubic potential
$V(\lambda;a,b)$ admits the \textit{simultaneous solution of two different quantization conditions},
namely $\sigma_{\pm 2}=0$.
\end{Cor*}

\subsection{Asymptotic values}\label{section:values}

As it was previously observed, Stokes multipliers are defined by particular normalized solutions
of equations (\ref{sys:1}) and (\ref{eq:schr}).
Following Nevanlinna, we define the monodromy data of equation (\ref{eq:schr}) in a more invariant way.

\begin{Def}
Let $\left\lbrace \varphi,\chi \right\rbrace$ be a basis of solution of (\ref{eq:schr}).

We  call \begin{equation}\label{def:wk}
 w_k(\varphi,\chi)=\lim_{\substack{\lambda\ \to \infty \\
\left|\arg{\lambda}-\frac{2 \pi k}{5}\right| < \frac{\pi}{5} -\varepsilon}}\frac{\varphi(\lambda)}{\chi(\lambda)} \in \mathbb{C}
\cup \infty \, , \; k \in \mathbb{Z}_5 \; .
\end{equation}
the k-th asymptotic value.
\end{Def}

We collect the main properties of the asymptotic values in the following
\begin{Lem}

\begin{itemize}
 \item[(i)]  Let $\varphi'= a \, \varphi + b \, \chi$ and $\chi'=c \, \varphi + d \, \chi'$,
$a,b,c,d \in \mathbb{C}$. Then
\begin{equation}\label{eq:moebius}
w_k(\varphi',\chi')= \frac{a \, w_k(\varphi,\chi) + b}{c \, w_k(\varphi,\chi)+d} \; .
\end{equation}
\item[(ii)] $w_{k-1}(\varphi,\chi)=w_{k+1}(\varphi,\chi)$ \textit{iff} $\sigma_k=0$ .
\item[(iii)]$w_{k+1}(\varphi,\chi) \neq w_{k}(\varphi,\chi)$

\end{itemize}

\begin{proof}
See \cite{elfving}.
\end{proof}
\end{Lem} 
Making use of equation (\ref{def:wk}), given the Stokes multipliers it is possible to calculate the
asymptotic values. The converse is also true.
In particular, the asymptotic values associated to the tritronqu\'ee int\'egrale can be chosen to be

\begin{equation}\label{tritronvalue}
w_0=0,w_1=w_{-2}=1,w_2=w_{-1}=\infty \; .
\end{equation}

\section{Stokes Complexes}\label{section:complexes}

In the complex WKB method a prominent role is played by the \textit{Stokes and anti-Stokes lines},
and in particular by the topology of the \textit{Stokes complex}, which is the union of
the Stokes lines.

The main result of this section is the Classification Theorem, where
we show that the topological classification of Stokes complexes
divides the space of cubic potentials into seven disjoint subsets.

Even though Stokes and anti-Stokes lines are well-known objects, there is no standard convention about
their definitions, so that
some authors call Stokes lines what others call anti-Stokes lines. 
We follow here the notation of Fedoryuk
\cite{fedoryuk}.

\begin{Rem*}
To simplify the notation and avoid repetitions, we study the Stokes lines only.
Every single statement in the following section remains true if
the word Stokes is replaced with the word anti-Stokes, provided
in equation (\ref{def:stokesline}) the angles $\varphi_k$ are replaced
with the angles $\varphi_k + \frac{\pi}{5}$.
\end{Rem*}

\begin{Def}\label{def:preliminar}
A simple (resp. double, resp. triple) zero $\lambda_i$ of $V(\lambda)=V(\lambda;a,b)$ is called a
simple (resp. double, resp. triple) turning point. All other points are called generic.

Fix a generic point $\lambda_0$ and a choice
of the sign of $\sqrt{V(\lambda_0)}$. We call action the analytic function
$$S(\lambda_0,\lambda)=\int_{\lambda_0}^{\lambda}\sqrt{V(u)}du$$ defined on 
the universal covering of $\lambda$-plane minus the turning points.
\end{Def}

Let $\tilde{i}_{\lambda_0}$  be the level curve of
the real part of the action
passing through a lift of $\lambda_0$.
Call its projection to the punctured plane $i_{\lambda_0}$.
Since $i_{\lambda_0}$ is
a one dimensional manifold, it is diffeomorphic to a circle or to a line.
If $i_{\lambda_0}$ is diffeomorphic to the real line, we choose one diffeomorphism
$i_{\lambda_0}(x),  x \in \mathbb{R}$
in such a way that the continuation along the curve of the
imaginary part of the action is a monotone increasing function of
$x \in \mathbb{R}$.

\begin{Lem}\label{lem:extension}
Let $\lambda_0$ be a generic point. Then $i_{\lambda_0}$ 
is diffeomorphic to the real line,
the limit $\lim_{x \to + \infty}i_{\lambda_0}(x)$ exists (as a point in $\mathbb{C} \bigcup \infty$)
and it satisfies the following dichotomy:

\begin{itemize}
 \item[(i)]Either $\lim_{x \to +\infty}i_{\lambda_0}(x)= \infty$
and the curve is asymptotic to one of the following rays of the complex plane
\begin{equation}\label{def:stokesline}
\lambda=\rho e^{i \varphi_k}, \varphi_k = \frac{(2k+1) \pi}{5} \,
, \rho \in \mathbb{R}^+, k \in \mathbb{Z}_5 \; ,
\end{equation}

\item[(ii)]or  $\lim_{x \to +\infty} i_{\lambda_0}(x)=\lambda_i$, where $\lambda_i$ is a turning point.
\end{itemize}

Furthermore,

\begin{itemize}
 \item[(iii)] if $\lim_{x \to \pm \infty}i_{\lambda_0}(x)=\infty$ then
the asymptotic ray in the positive direction is different from the asymptotic ray in
the negative direction.

\item[(iv)] Let $\varphi_k, k \in \mathbb{Z}_5$ be defined as in equation (\ref{def:stokesline}).
Then $\forall \varepsilon >0,\exists K \in \mathbb{R}^+ $
such that if $\varphi_{k-1} + \varepsilon<\arg{\lambda_0}<\varphi_k - \varepsilon
\mbox{ and } \av{\lambda_0}>K  $,
then $\lim_{x \to \pm \infty}i_{\lambda_0}(x)=\infty$. Moreover
the asymptotic rays of $i_{\lambda_0}$ are the ones with arguments $\varphi_{k}$
and $\varphi_{k-1}$.
\end{itemize}

\begin{proof} See \cite{strebel}.
 
\end{proof}

\end{Lem}

\begin{Def}
We call Stokes line the trajectory of any curve $i_{\lambda_0}$  such
that there exists at least one turning point belonging to its boundary.

We call a Stokes line internal if
$\infty$ does not belong to its boundary.

We call Stokes  complex the union of all the Stokes lines
together with the turning points.
\end{Def}

We state all important properties of the Stokes lines in the
following

\begin{Thm}\label{thm:curves}
The following statements hold true 
\begin{itemize}

 \item[(i)] The Stokes complex
 is simply connected. In particular, the boundary of any internal Stokes  line
 is the union of two different turning points.

\item[(ii)] Any simple (resp. double, resp. triple) turning point belongs to the boundary
of 3 (resp. 4, resp 5) Stokes lines.

\item[(iii)]If a turning point belongs to the boundary of two different non-internal Stokes
lines then these lines have different asymptotic rays.

\item[(iv)]For any ray with the argument  $\varphi_k$ as in equation (\ref{def:stokesline}), there
exists a Stokes line asymptotic to it.

\end{itemize}

\begin{proof}See \cite{strebel}.
\end{proof}

\end{Thm}

\subsection{Topology of Stokes complexes}

In what follows, we give a complete classification of the Stokes complexes, with 
respect to the orientation preserving homeomorphisms of the plane.

We define the map $L$ from the $\lambda$-plane to the interior of the unit disc as
\begin{eqnarray}L &:& \mathbb{C} \to D_1 \nonumber \\
 L(\rho e^{i \varphi}) &=& \frac{2}{\pi} e^{i \varphi}\arctan{\rho}.
\end{eqnarray}

The image under the map $L$ of the Stokes complex is naturally a decorated graph
embedded in the closed unit disc. The vertices are
the images of the turning points and the five points
on the boundary of the unit disc with arguments $\varphi_k$, with $\varphi_k$
as in equation (\ref{def:stokesline}). The bonds are obviously the images
of the Stokes lines.
We call
the first set of vertices \textit{internal} and the second set of vertices \textit{external}.
External vertices are decorated with the numbers $k \in \mathbb{Z}_5$.
We denote $\cal{S}$ the decorated embedded graph just described.
Notice that due to Theorem \ref{thm:curves} (iii), there exists not more
than one bond connecting two vertices.

The combinatorial properties of $\cal{S}$ are described in the following

\begin{Lem}\label{lem:graph}
$\cal{S}$ possesses the following properties
\begin{itemize}
 \item[(i)] the sub-graph spanned by the internal vertices has no cycles.
\item[(ii)] Any simple (resp. double, resp. triple) turning point has valency
3 (resp. 4, resp. 5).
\item[(iii)] The valency of any external vertex is at least one. 
\end{itemize}
\begin{proof}
\begin{itemize}
 \item[(i)] Theorem \ref{thm:curves} part (i)
\item[(ii)] Theorem \ref{thm:curves} part (ii)
\item[(iii)] Theorem \ref{thm:curves} part (iv)
\end{itemize}

\end{proof}

\end{Lem}

\begin{Def}
We call an {\rm admissible graph} any decorated simple graph embedded in the closure of the unit disc,
with three internal vertices and five decorated external vertices, such that
(i) the cyclic-order inherited from the decoration coincides
with the one inherited from the counter-clockwise orientation of
the boundary, and (ii) it satisfies all the
properties of Lemma \ref{lem:graph}.
We call two admissible graphs equivalent if there exists an orientation-preserving homeomorphism
of the disk mapping one graph into the other.
\end{Def}

\begin{figure}

\input{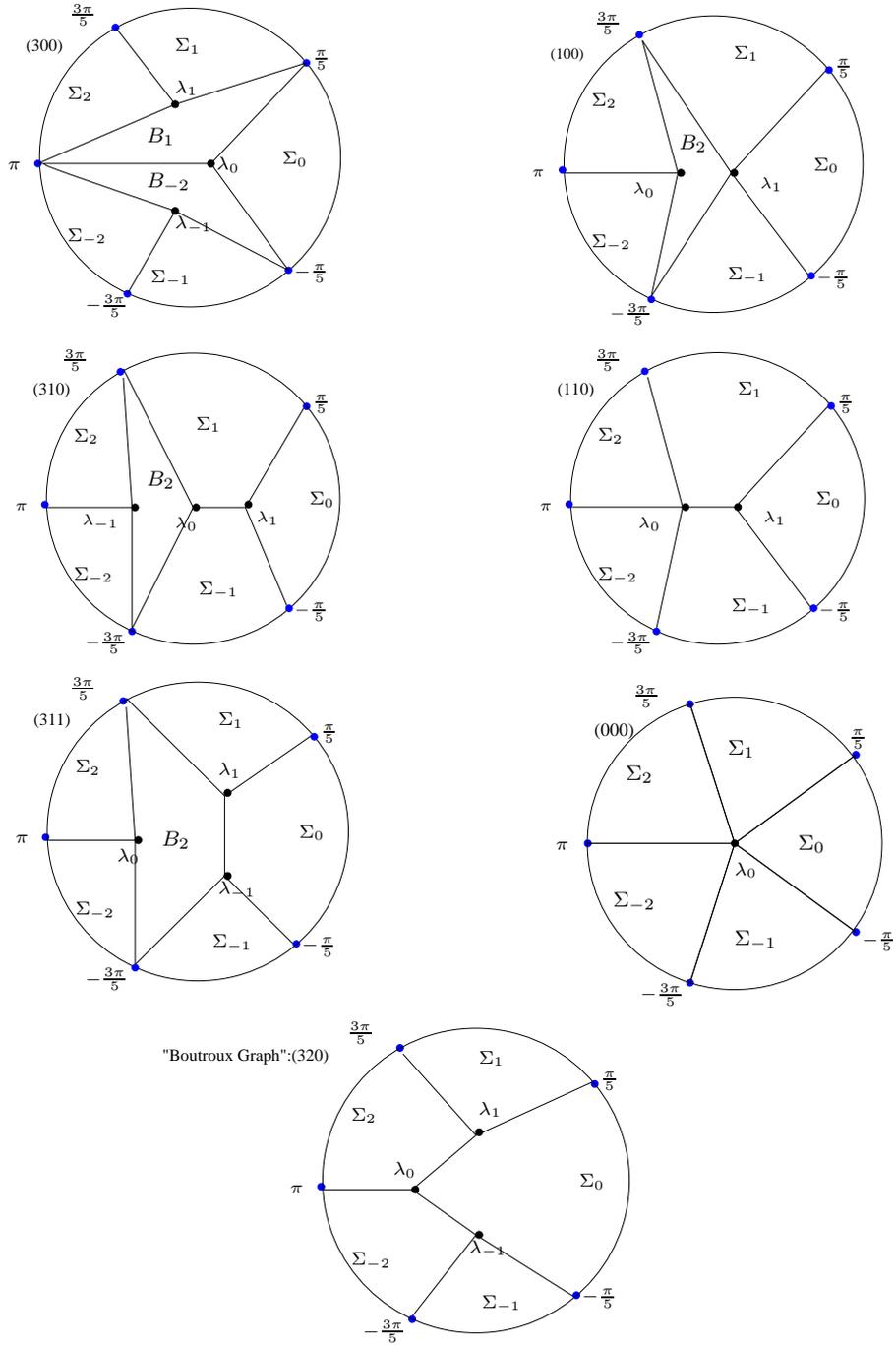}
 
\caption{All the equivalence classes of admissible graphs.}
\label{figure:maps}

\end{figure}

\begin{Thm}\label{thm:classification}Classification Theorem

All equivalence classes of admissible maps are,
modulo a shift $k \to k +m, m \in \mathbb{Z}_5$ of the decoration,
the ones depicted in Figure \ref{figure:maps}.

\begin{proof}

Let us start analyzing the admissible graphs with three internal vertices and no internal edges.

Any internal vertex is adjacent to a triplet of external vertices.
Due to the Jordan curve theorem,
there exists an internal vertex, say $\lambda_0$,
adjacent to a triplet of non consecutive external vertices.
Performing a shift, they can be chosen to be the ones labelled by $0, 2, -1$.
Call the respective edges $e_0,e_{-1},e_2$.

The disk is cut in three disjoint domains by those three edges. No internal vertices
can belong to the domain cut by 
$e_0$ and $e_4$, since it could be adjacent only to
two external vertices,  namely the ones labelled with $0$ and $-1$.
By similar reasoning it is easy to show that one and only one vertex belong
to each remaining domains.

Such embedded graph is equivalent to the graph (300).

Classifications for all other cases may proved by similar methods.

\end{proof}

\end{Thm}

The equivalence classes are encoded by a triplet of numbers (a b c):
$a$ is the number of simple turning  points,
$b$ is the number of internal Stokes  lines, while $c$ is a progressive number,
distinguishing non-equivalent graphs with same $a$ and $b$. Some additional
information shown Figure \ref{figure:maps} will be explained in the next section.

\begin{Rem*}
For any admissible graph there exists a real polynomial  with
an equivalent Stokes complex.
\end{Rem*} 
\begin{Rem*}
Notice that the automorphism group of every graph in Figure \ref{figure:maps}
is trivial.
Therefore the unlabelled vertices can be labelled.
In the following
we will label the turning points as in figure \ref{figure:maps}.
We denote "Boutroux graph" the graph $(320)$
\end{Rem*}

\subsection{Stokes Sectors}

In the $\lambda$-plane the complement of the Stokes complex is the disjoint union of a 
finite number of connected and simply-connected domains, each of them called a \textit{sector}.

Combining Theorem \ref{thm:curves} and the Classification Theorem we
obtain the following
\begin{Lem}
All the curves $i_{\lambda_0}$, with $\lambda_0$ belonging to a given sector, have the same two asymptotic rays.
Moreover, two different sectors have different pairs of asymptotic rays.

For any $k \in \mathbb{Z}_5$ there is a sector, called the k-th \textit{Stokes sectors}, whose
asymptotic rays have arguments $\varphi_{k-1}$ and $\varphi_k$. 
This sector will be denoted $\Sigma_k$. The boundary  $\partial \Sigma_k$ of each $\Sigma_k$
is connected.

Any other sector has asymptotic rays with arguments $\varphi_{k-1}$ and $\varphi_{k+1}$, for some $k$.
We call such a sector the k-th sector of band type, and we
denote it $B_k$. The boundary $\partial B_k$ of each $B_k$ has two connected components.

\end{Lem}

Choose a sector and a point $\lambda_0$ belonging to it.
The function $S(\lambda_0,\lambda)$ is easily seen to be bi-holomorphic into the image
of this sector.
In particular, with one choice of the sign of $\sqrt{V}$
it maps a Stokes sector into the half plane ${\rm Re}S > c$, for some $-\infty<c<0$
while it maps a $B_k$ sector in the vertical strip $ c<{\rm Re}S <d$,
for some $-\infty<c<0<d<+\infty$.

\begin{Def}
We call
a differentiable curve $\gamma:[0,1] \to \mathbb{C}$ an admissible path
provided $\gamma$ is injective on $[0,1[$,
$\lambda_i \notin \gamma([0,1])$, for all turning points  $\lambda_i$,
and ${\rm Re}S(\gamma(0), \gamma(t))$ is a monotone function of $t \in [0,1]$.

We say that $\Sigma_j \leftrightarrows \Sigma_k$ if there exist $\mu_j \in \Sigma_j$,
$\mu_k \in \Sigma_k$ and an admissible path such that $\gamma(0)=\mu_j, \gamma(1)=\mu_k$.

\end{Def}

The relation $\leftrightarrows$ is obviously reflexive and symmetric but it is not in general transitive.

Notice that $\Sigma_j \leftrightarrows \Sigma_k$ \textit{if and only if} for every point $\mu_j \in \Sigma_j$
and every point $\mu_k \in \Sigma_k$ an admissible path exists.

\begin{Lem}\label{lem:relation}
The relation $\leftrightarrows$ depends only on the equivalence class of the
Stokes complex $\cal{S}$.

\begin{proof}
Consider an admissible path from $\Sigma_j$ to $\Sigma_k$, $j \neq k$.
The path is naturally associated to the sequence of Stokes lines that it crosses.
We denote the sequence $l_n, n=0,\dots, N $, for some $N \in \mathbb{N}$.
We continue analytically $S(\mu_j, \cdot)$ to a covering of the union of the Stokes sectors
crossed by the path together with the Stokes lines belonging
to the sequence.
Since $S(\mu_j, \cdot)$ is constant along each 
connected component of the boundary of every lift of a sector crossed by the path,
then each of such connected components cannot be crossed twice by the path.
Hence, due to the classification theorem no admissible path is a loop.
Therefore, the union of the Stokes sectors
crossed by the path together with the Stokes lines belonging
to the sequence is simply connected.

Conversely, given any injective sequence of Stokes lines $l_n, \, n=0\dots,N$
such that
for any $0 \leq n \leq N-1$, $l_n$ and $l_{n+1}$ belong to two different connected
components of the boundary of a same sector, there
exists an admissible path with that associated sequence.
This last observation implies that the relation $\leftrightarrows$ depends only
on the topology of the graph $\cal{S}$. Moreover, if the sequence exists it is unique;
indeed, if there existed two admissible paths, joining the same $\mu_j$ and $\mu_k$
but with different sequences, then there would be  an admissible loop.
\end{proof}
\end{Lem}

\begin{table}
\begin{tabular}{|c|c|}
\hline Map & Pairs of non consecutive Sectors not satisfying the relation $\leftrightarrows$ \\ 
\hline 300 & None \\ 
\hline 310 & $\left(\Sigma_0, \Sigma_{2}\right)$, $\left( \Sigma_0, \Sigma_{-2} \right)$  \\ 
\hline 311 & $\left(\Sigma_1, \Sigma_{-1}\right)$ \\ 
\hline 320 & $\left(\Sigma_{1}, \Sigma_{-1}\right)$, $\left(\Sigma_{1}, \Sigma_{-2}\right)$, $\left(\Sigma_{-1}, \Sigma_{2}\right)$ \\ 
\hline 100 & $\left(\Sigma_{1}, \Sigma_{-1}\right)$, $\left(\Sigma_{0}, \Sigma_{-2}\right)$, $\left(\Sigma_{0}, \Sigma_{2}\right)$ \\ 
\hline 110 & All but $\left(\Sigma_{1}, \Sigma_{-1}\right)$ \\ 
\hline 000 & All \\ 
\hline  
\end{tabular}
\caption{Computation of the relation $\leftrightarrows$}
\label{table:prec}
\end{table} 

With the help of Lemma \ref{lem:relation} and of the Classification Theorem,
relation $\leftrightarrows$ can be easily computed, as it is
shown in Table \ref{table:prec}. As it is evident from Figure \ref{figure:maps},
for any graph type we have that $\Sigma_k \leftrightarrows \Sigma_{k+1},
\, \forall k \in \mathbb{Z}_5$.

\section{Complex WKB Method and Asymptotic Values}\label{section:wkb}
In this section we introduce the \textit{WKB functions} $j_k, k \in \mathbb{Z}_5$ and use them
to evaluate the asymptotic values of
equation (\ref{eq:schr}). The topology of the Stokes complex will
show all its importance in these computations.

On any Stokes sector $\Sigma_k$, we define the functions

\begin{eqnarray} \label{def:Sk}
S_k(\lambda)&=&S(\lambda^*,\lambda) \; , \\ \label{def:Lk}
L_k(\lambda)&=&-\frac{1}{4}\int_{\lambda^*}^{\lambda}\frac{V'(u)}{V(u)}du \; ,\\ \label{def:jk}
j_k(\lambda)&=&e^{-S_k(\lambda)+L_k(\lambda)} \; .
\end{eqnarray}

Here $\lambda^*$ is an arbitrary point belonging to $\Sigma_k$ and
the branch of $\sqrt{V}$ is such that ${\rm Re}S_k(\lambda)$ is bounded from below.

We call $j_k$ the k-th \textit{WKB function}.

\subsection{Maximal Domains}\label{sub:domains}

In this subsection we construct the k-th \textit{maximal domain}, that
we denote  $D_k$. This is
the domain of the complex plane where the k-th WKB function approximates
a solution of equation (\ref{eq:schr}).

\begin{figure}[htbp]
\begin{center}
\input{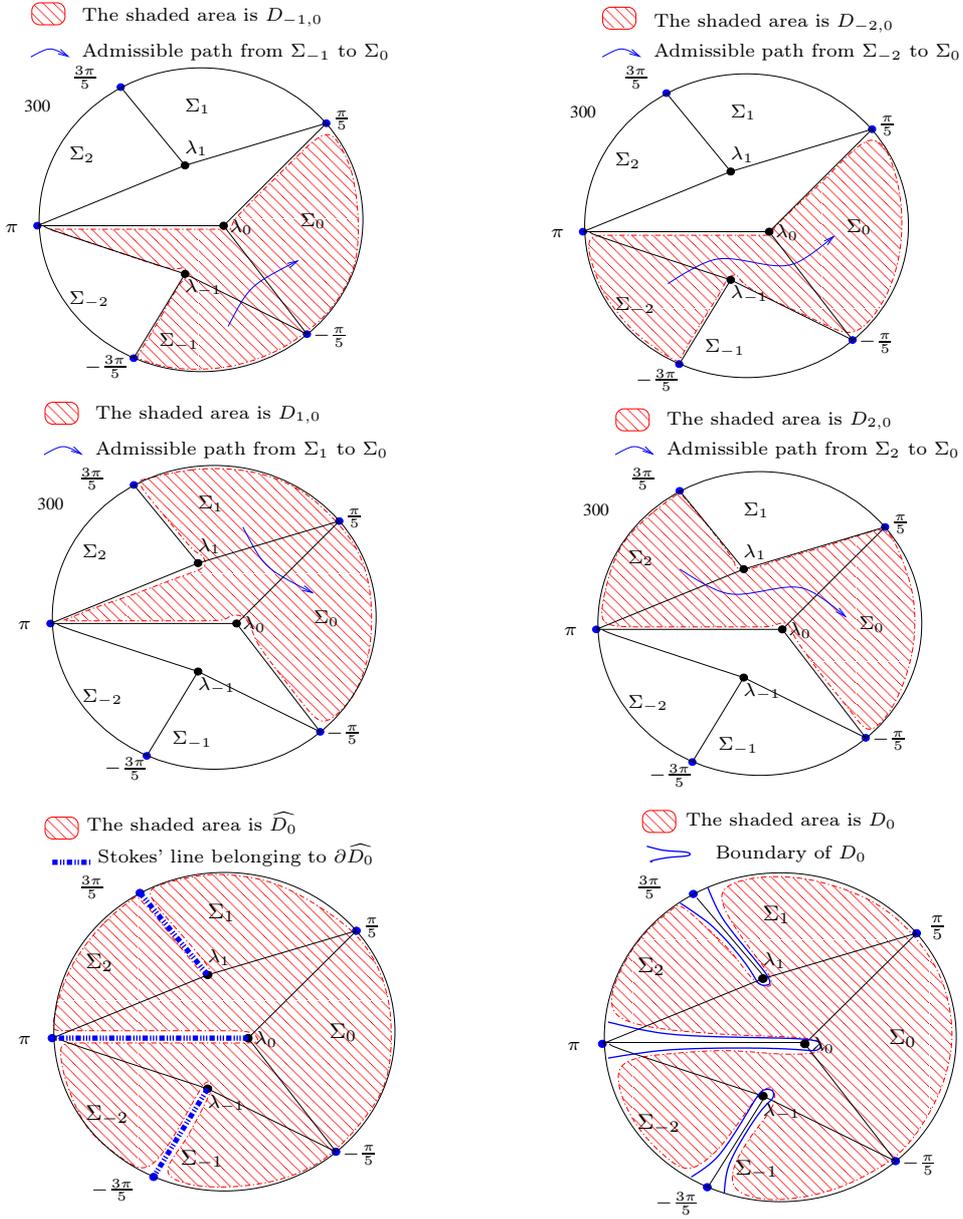}

\end{center}

\caption{In the drawings, the construction of $D_0$ for a graph of type $(300)$ is depicted.}
\label{figure:Dk}
\end{figure}

The construction is done for any $k$ in a few steps (see Figure \ref{figure:Dk} for the example
of the Stokes complex of type (300)):

\begin{itemize}
\item[(i)] for every $\Sigma_l$ such that $\Sigma_l \leftrightarrows \Sigma_k$, denote $D_{k,l}$
the union of the sectors and of the Stokes lines crossed by any admissible path connecting $\Sigma_l$
and $\Sigma_k$.

\item[(ii)]Let $\widehat{D_k}= \bigcup_{l}D_{l,k}$. Hence $\widehat{D_k}$ is a connected and simply
connected 
subset of the complex plane whose boundary $\partial \widehat{D_k}$ is the union of some Stokes lines. 

\item[(iii)]Remove a $\delta \mbox{-tubular}$ neighborhood of the boundary $\partial \widehat{D_k}$, for an arbitrarily small $\delta>0$, such that the resulting domain is still connected.

\item[(iii)]For all $l\neq k, l \neq k-1$,
remove from $\widehat{D_k}$ an angle
$\lambda= \rho e^{i \varphi}, \left| \varphi-\varphi_l \right| < \epsilon, \rho>R$,
for $\e$ arbitrarily small and $R$ arbitrarily big, in such a way that the resulting domain
is still connected. The remaining domain is $D_k$.

\end{itemize}

\subsection{Main Theorem of WKB Approximation}

We can now state the main theorem
of the WKB approximation.

\begin{Thm}[G.D. Birkhoff \cite{birkhoff}, Olver \cite{olver}]\label{thm:birkh}
Continue the WKB function $j_k$ to $D_k$.
Then there exists a solution $\psi_k(\lambda)$ of (\ref{eq:schr}), such that for all $\lambda \in D_k$
\begin{eqnarray*}
\left| \frac{\psi_k(\lambda)}{j_k(\lambda)} -1 \right| & \leq& g(\lambda) \left( e^{2\rho(\lambda)}-1 \right)\\
\left| \frac{\psi_k'(\lambda)}{j_k(\lambda)\sqrt{V(\lambda)}} +1 \right|
&\leq&  \left|\frac{V'(\lambda)}{4V(\lambda)^{\frac{3}{2}}}\right|+
(1+\left|\frac{V'(\lambda)}{4V(\lambda)^{\frac{3}{2}}}\right|)g(\lambda)( e^{2 \rho(\lambda)} -1) 
\end{eqnarray*}

Here $\rho_k$ is a bounded positive continuous function, called the \textit{error function},
satisfying

$$\lim_{\substack{\lambda \to \infty \\
\varphi_{k-1}<arg\lambda < \varphi_{k+1} }} \rho_k(\lambda) =0 \;, $$

and $g(\lambda)$ is a positive function such that $g(\lambda) \leq 1$ and

$$\lim_{\substack{\lambda \to \infty \\
\lambda \in D_k \cap \Sigma_{k \pm 2}}} g(\lambda)=\frac{1}{2} \; . $$

\begin{proof}
 The proof is in the appendix \ref{app:A}.
\end{proof}

\end{Thm}

Notice that $j_k$ is \textit{sub-dominant} (i.e. it decays exponentially) in $\Sigma_k$ 
and \textit{dominant} (i.e. it grows exponentially) in $\Sigma_l, \forall l \neq k$.

For the properties of the error function, $\psi_k$ is subdominant in $\Sigma_k$
and dominant in $\Sigma_{k \pm1}$. Therefore, in any Stokes sector $\Sigma_k$
there exists a subdominant solution, which is defined uniquely up to multiplication
by a non zero constant.

\subsection{Computations of Asymptotic Values in WKB Approximation} \label{par:values}

The aim of this paragraph is to compute the asymptotic values for
the Schr\"odinger equation (\ref{eq:schr}) in WKB approximation.
We explicitely work out the example of the Stokes complex of type $(320)$,
relevant to the study of poles of the int\'egrale tritronqu\'ee.

\begin{Def}\label{Def:errors}
Define the relative errors
\begin{equation*}
\rho_l^k = \left\lbrace 
             \begin{aligned}
               \lim_{\substack{\lambda \to \infty \\ \lambda \in
\Sigma_{k} \cap D_{l}}}\rho_l(\lambda) , \; \mbox{ if } \, \Sigma_l \leftrightarrows \Sigma_k  \\
 \infty , \qquad \quad \; \mbox{otherwise} 
             \end{aligned}
\right. 
\end{equation*}
and the asymptotic values
\begin{equation} \label{Def:values}
w_k(l,m) \stackrel{def}{=} w_k(\psi_l,\psi_m).
\end{equation}
We say that $\Sigma_k \sim \Sigma_l$ provided $\rho_l^k < \frac{\log 3}{2}$.
The relation $\sim$ is a sub-relation of $\leftrightarrows$.
\end{Def}
Notice that $\rho_l^{l+1}=0$  and  $\rho_l^m=\rho_m^l$ (see Appendix \ref{app:A}).

In order to compute the asymptotic value $w_k(l,m)$,
we have to know the asymptotic behavior
of $\psi_l$ and $\psi_m$ in $\Sigma_k$. By Theorem \ref{thm:birkh},
$$\lim_{\substack{\lambda \to \infty \\ \lambda \in
\Sigma_{k} \cap D_{l}}}\frac{\psi_l(\lambda)}{j_l(\lambda)}
\neq 0 \, , \mbox{ if } \frac{1}{2}(e^{2 \rho_k^l}-1) < 1 \; . $$
Hence the asymptotic behavior of $\psi_l$  in $\Sigma_k$ can be related
to the asymptotic behavior of $j_l$ in $\Sigma_k$  if the relative error $\rho_l^k$ is so small that
the above inequality holds true, i.e. if $\Sigma_k \sim \Sigma_l$.

\begin{Rem*}
Depending on the type of the graph $\cal{S}$, there may not exist
two indices $k \neq l$ such that all the relative errors
$\rho_l^n,\rho_k^n, n \in \mathbb{Z}_5$
are small. However it is often possible to compute
an approximation of all the asymptotic values $w_n(l,k)$ using the strategy below.
\end{Rem*}

\begin{itemize}
 \item[(i)] We select a pair of non consecutive Stokes sectors $\Sigma_l$, $\Sigma_{l+2}$, 
 with the hypothesis that the functions $\psi_l$ and $\psi_{l+2}$ are linearly independent, so that
$w_l(l,l+2)=0, w_{l+2}(l,l+2)=\infty$. 
Since $\rho_l^{l+1}=\rho_{l+2}^{l+1}=0$ then
$$w_{l+1}(l,l+2)=\lim_{\substack{\lambda \to \infty \\ \lambda \in
\Sigma_{l+1} \cap D_{l} \cap D_{l+2} } } \frac{j_l(\lambda)}{j_{l+2}(\lambda)} \; .$$
Therefore, we find three exact and distinct asymptotic values.

\item[(ii)] For any $k \neq l,l+1,l+2$ such that $\Sigma_l \sim \Sigma_k$ and $\Sigma_{l+2} \sim \Sigma_k$,
we define the approximate asymptotic value
$$
\widehat{w}_k(l,m)=\lim_{\substack{\lambda \to \infty \\ \lambda \in
\Sigma_{k} \cap D_{l} \cap D_{l+2} }}\frac{j_l(z)}{j_m(z)} \; .
$$

The spherical distance between $w_k(l,m)$ and $\widehat{w}_k(l,m)$ may be easily estimated
from above knowing the relative errors $\rho_k^l$ and $\rho_k^{l+2}$.


If for any $k \neq l,l+1,l+2$,  $\Sigma_l \sim \Sigma_k$ and $\Sigma_{l+2} \sim \Sigma_k$, then
the calculation is completed.

\item[(iii)]If not, we can use the fact that  quintuplets of asymptotic values
for different choices of $\psi_l, \psi_{l+2}$ are related by a
M\"obius transformation (see formula (\ref{eq:moebius})).
If for some pair $(l,l+2)$ the assumption $\Sigma_l \sim \Sigma_k, \,\Sigma_{l+1} \sim \Sigma_k$
fails to be true for just one value of the index $k=k^*$, and, for another pair $(l',l'+2)$
the assumption $\Sigma_{l'} \sim \Sigma_{k'} ,\, \Sigma_{l'+2} \sim \Sigma_{k'}$
fails to be true for just one valued of the index $k'=k'^*$, with $k'^* \neq k^*$,
then there are three values of the index $m \in \mathbb{Z}_5$
such that an approximation of $w_m(l,l+2)$ and $w_m(l',l'+2)$ is computable.
Since any M\"obius transformation is fixed by the action on three values, then
we can compute an approximation of the transformation relating the quintuplets
$w_k(l,l+2)$ and $w_k(l',l'+2)$ for any $k \in \mathbb{Z}_5$.
Hence, we can calculate an approximation of the whole
quintuplets $w_k(l,l+2)$ and $w_{k'}(l',l'+2)$.

\end{itemize}

\begin{Rem*}
As shown in Table \ref{table:prec}, the relation $\leftrightarrows$ is uniquely characterized by the
graph type.
For the sake of computing the asymptotic values the important relation
is $\sim$ and not $\leftrightarrows$.
Indeed, the calculations for a given graph type, say $(a\; b \; c)$,
are valid for (and only for) all the potentials whose relation $\sim$
is equivalent to the relation $\leftrightarrows$ characterizing the graph type $(a \; b\; c)$.
\end{Rem*}

Due to the above remark, in what follows we suppose that the relation
$\sim$ is equivalent to the relation $\leftrightarrows$.
We have the following
\begin{Lem}\label{lem:exclusion}
Let $V(\lambda;a,b)$ such that the type of the Stokes complex is $(300)$, $(310)$, $(311)$; moreover,
suppose that the $\sim$ relation coincides with $\leftrightarrows$. Then
all the asymptotic values  of equation (\ref{eq:schr})  are pairwise distinct,
but for at most pair.
\begin{proof}
For a graph of type $(300)$ or $(311)$ the thesis is trivial.
For a graph of type $(320)$, it may be that $w_0=w_2$ or $w_0=w_{-2}$. Since
$w_2 \neq w_{-2}$ the thesis follows.
\end{proof}

\end{Lem}

We completely work out the case of Stokes complex of type $(320)$, while
for the other cases we present the results only. Due to Lemma \ref{lem:exclusion},
we omit the results for potentials whose graph type is $(300)$, $(310)$ and $(311)$.

\paragraph{Boutroux Graph = 320}

We suppose that $\Sigma_0 \sim \Sigma_{\pm2}$.

\begin{figure}[htbp]
\begin{center}
\input{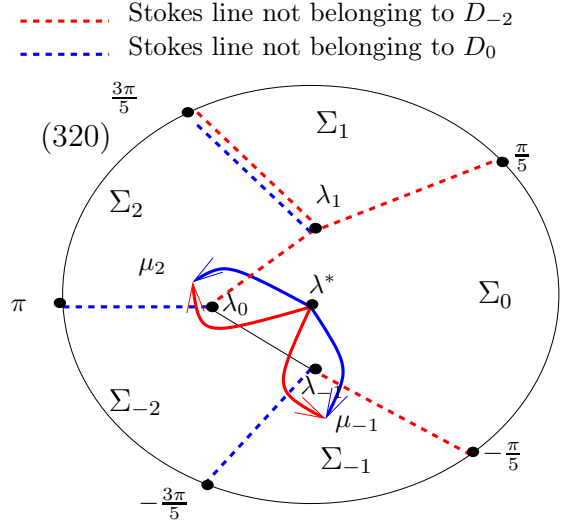}
 \end{center}

\caption{Calculation of $w_{-1}(0,-2)$ and of $\widehat{w}_2(0,-2)$}
\label{figure:320calc}
\end{figure}

Let us consider first the pair $\Sigma_0$ and $\Sigma_{-2}$.
In Figure \ref{figure:320calc} the maximal domains $D_0$ ans $D_{-2}$
are depicted by colouring the Stokes lines
not belonging to them blue and red respecitvely.
In particular
$S_0, L_0
,j_0$ (resp. $S_{-2},L_{-2},j_{-2}$)
can be extended to all $D_0$ (resp. $D_{-2}$) along any curve
that does not intersect any blue (resp. red) Stokes line.

We fix a point $\lambda^* \in \Sigma_0$ such
that $S_0(\lambda^*)=S_{-2}(\lambda^*)=L_0(\lambda^*)=L_{-2}(\lambda^*)=0$.

By definition 
\begin{eqnarray*}
 \widehat{w}_k(0,-2) &=& \lim_{\lambda \to \infty_k}\frac{j_0(\lambda)}{j_{-2}(\lambda)}  \\
 &=&\lim_{\lambda \to \infty_k}e^{-S_0(\lambda)+S_{-2}(\lambda)} e^{L_0(\lambda)-L_{-2}(\lambda)} \quad ,
\end{eqnarray*}
Here $\lambda \to \infty_k$ is a short-hand notation for $\lambda \to \infty ,
\lambda \in \Sigma_k \cap D_0 \cap D_{-2}$.
We calculate $\widehat{w}_k(0,-2)$ for $k=-1,2$.

We first calculate $\lim_{\lambda \to \infty_k}e^{-S_0(\lambda)+S_{-2}(\lambda)}$.

Notice that
$ \frac{\partial S_0}{\partial \lambda } =\frac{\partial S_{-2}}{\partial \lambda }$
in $\Sigma_k$.
Hence
$$\lim_{\lambda \to \infty_k}-S_0(\lambda)+S_{-2}(\lambda)= -S_0(\mu_k)+S_{-2}(\mu_k) \, , k=-1,2 \; ,$$
where $\mu_k$ is any point belonging to $\Sigma_k$
(in Figure \ref{figure:320calc}, the paths of integration defining $S_0(\mu_k)$
and $S_{-2}(\mu_k)$ are coloured  blue and red respectively).

On the other hand, since $ \frac{\partial S_0}{\partial \lambda } =-\frac{\partial S_{-2}}{\partial \lambda }$
in $\Sigma_0 \bigcup \Sigma_{-2}$, we have that
$$ -S_0(\mu_k)+S_{-2}(\mu_k) = -2 S_0(\lambda_s) \, , s=-1 \mbox{ if } k=-1
\mbox{ and } s=0 \mbox{ if } k=2 \; . $$ 

We now compute $\lim_{\lambda \to \infty_k}e^{L_0(\lambda)-L_{-2}(\lambda)}$.
Since $ \frac{\partial L_0}{\partial \lambda } =\frac{\partial L_{-2}}{\partial \lambda }$
in $D_0 \bigcap D_{-2} $, we have that
\begin{eqnarray*}
\lim_{\lambda \to \infty_k}L_0(\lambda)-L_{-2}(\lambda) &=& L_0(\mu_k)-L_{-2}(\mu_k) \,\; , \\
L_0(\mu_k)-L_{-2}(\mu_k) &=& -\frac{1}{4}\oint_{c_k}\frac{V'(\mu)}{V(\mu)}d\!\mu \, , k=-1,2 \; .
\end{eqnarray*}
Here $c_k$ is the blue path connecting $\lambda^*$ with $\mu_k$ composed with the inverse
of the red path connecting $\lambda^*$ with $\mu_k$ (see Figure \ref{figure:320calc}).

Therefore, we have
$$ \lim_{\lambda \to \infty_k}L_0(\lambda)-L_{-2}(\lambda)= -\sigma \frac{2 \pi i}{4}
\, , \sigma=-1 \mbox{ if } k=-1 \mbox{ and } \sigma=+1 \mbox{ if } k=2 \; .  $$

Combining the above computations, we get

\begin{equation*}
 w_{-1}(0,-2)= i \, e^{-2 S_0(\lambda_{-1})} , \; \widehat{w}_2(0,-2)= - i\, e^{-2 S_0(\lambda_0)} .
\end{equation*}
We stress that $w_{-1}(0,-2)$ is exact while $\widehat{w}_2(0,-2)$ is an approximation.

Performing the same computations for the pair $\Sigma_0$ and $\Sigma_{2}$, we obtain

\begin{equation*}
 w_{1}(0,2)=- i \, e^{-2 S_0(\lambda_{1})} , \; \widehat{w}_{-2}(0,2)=  i\, e^{-2 S_0(\lambda_0)} 
\end{equation*}

Having calculated the triplet of asymptotic values $w_0,w_2,w_{-2}$ for two different pairs of Stokes
sectors, we can compute an approximation of the M\"obius transformation relating all the asymptotic values
for the two pairs: 
\begin{equation*}
 \widehat{w}_{k}(0,-2)=- i\, e^{-2 S_0(\lambda_{0})} \frac{\widehat{w}_{k}(0,2)}
{\widehat{w}_{k}(0,2) - i \, e^{-2 S_0(\lambda_0)}} \, , \; k \in \mathbb{Z}_5 \quad .
\end{equation*}

We eventually compute the last asymptotic value for the pair $\Sigma_0, \Sigma_{-2}$, that is

\begin{eqnarray*}
  \widehat{w}_{1}(0,-2) &=&-i \, \frac{e^{-2 S_0(\lambda_{1})}}
{1+e^{-2 (S_0(\lambda_{1})-S_0(\lambda_0))}} \quad .
\end{eqnarray*}

\subparagraph{Quantization Conditions}

The computations above provides us with the following quantization conditions:

\begin{eqnarray}
 \widehat{w}_{1} &=& w_{-2} \Leftrightarrow
e^{-2 (S_0(\lambda_{1})-S_0(\lambda_0))} = -1 \label{eq:bohr1} \\
\widehat{w}_{2} &=& w_{-1} \Leftrightarrow
e^{-2 (S_0(\lambda_{-1})-S_0(\lambda_0))} = -1 \label{eq:bohr2} \\
\widehat{w}_{1} &=& w_{-1} \Leftrightarrow
e^{-2 (S_0(\lambda_{1})-S_0(\lambda_{-1}))} =
 -1+ e^{-2 (S_0(\lambda_{1})-S_0(\lambda_{0}))} \label{eq:bohr3}
\end{eqnarray}

We notice that equation (\ref{eq:bohr3}) is incompatible both with
(\ref{eq:bohr1}) and  (\ref{eq:bohr2}). Equations (\ref{eq:bohr1}) and  (\ref{eq:bohr2}) are
Bohr-Sommerfeld quantizations.

As was shown in equation (\ref{tritronvalue}), the poles of the integr\`ale tritronqu\`ee
are related to the polynomials such that $w_1=w_{-2}$ and $w_{-1}=w_2$.
Since equations (\ref{eq:bohr1}) and (\ref{eq:bohr2}) can be simultaneously solved,
solutions of system (\ref{eq:bohr1},\ref{eq:bohr2}) describe, in WKB approximation,
polynomials related to the int\'egrale tritronqu\'ee.
System (\ref{eq:bohr1},\ref{eq:bohr2})  was found by
Boutroux in \cite{boutroux} (through a completely different analysis),
to characterize the asymptotic distribution of the poles of the integr\`ale tritronqu\`ee.
Therefore we call (\ref{eq:bohr1},\ref{eq:bohr2}) the \textit{Bohr-Sommerfeld-Boutroux system}.

Equation (\ref{eq:bohr3}) will not be studied in this paper, even though is quite remarkable.
Indeed, it describes the breaking of the PT symmetry
(see \cite{delabaere} and \cite{bender}).















\paragraph{Case $(100)$}
\begin{eqnarray*}
 w_{0}(1,-1) &=& -1 \\
\widehat{w}_{-2}(1,-1) &=& \widehat{w}_{2}(1,-1) =1
\end{eqnarray*}

Since $w_0 \neq \widehat{w}_{\pm2}$ and $w_2 \neq w_{-2}$,
if the error $\rho_1^{-2}$ or $\rho_{-1}^2$ is small enough, then all
the asymptotic values are pairwise distinct.

\paragraph{Case $(110)$}

\begin{eqnarray*}
 \widehat{w}_{-1}(1,-2) &=& 1 \\
 w_{2}(1,-2) &=& - 1
\end{eqnarray*}

In this case, it is impossible to calculate $w_0$ with the WKB method
that has been here developed. Hence it may be that
either $w_0=w_2$ or $w_0 = w_{-2}$.

Notice, however, that $(110)$ is the graph only of a very restricted class of potentials namely
$V(\lambda)= (\lambda+\lambda_0)^2(\lambda -2 \lambda_0)$, where $\lambda_0$ is real and positive.
Since the potential is real then $w_0 \neq w_{\pm 2} $.

\paragraph{Case $(000)$}
In this case, no asymptotic values can be calculated.
Notice, however, that  $V(\lambda)=\lambda^3$ is the only potential
with graph $(000)$.
For this potential the asymptotic values can be computed exactly, simply using symmetry considerations.
Indeed one can choose $w_k=e^{\frac{2 k \pi}{5}i} \, , \; k \in \mathbb{Z}_5$.

\subsection{Small Parameter}\label{sec:symanzik}

The WKB method normally applies to problem with an \textit{external small parameter}, usually denoted
$\hbar$ or $\e$. In the study of the distributions of poles of a given
solution $y$ of P-I there is no external small parameter and we have to explore
the whole space of cubic potentials. The aim of this section is to introduce an
\textit{internal small parameter} in the space of cubic potentials, that
greatly simplifies our study.

On the linear space of cubic potentials in canonical form
\begin{equation*}
V(\lambda;a,b)= 4\lambda^3-2 a \lambda -28 b,
\end{equation*}

we define the following action of the group $\mathbb{R}^+ \times \mathbb{Z}_{5}$
(similar to what is called Symanzik rescaling in \cite{simon})

\begin{eqnarray}\label{def:action}
(x,m)[V(\lambda;a, b)] = V(\lambda; \Omega^{2m} x^2 a,\Omega^{3m} x^3 b),
 \, x \in \mathbb{R}^+, \, m \in \mathbb{Z}_5,\, \Omega=e^{\frac{2 \pi}{5}i}.
\end{eqnarray}

The induced action on the graph $\cal{S}$, on the relative error $\rho_l^m$, and on
the difference $S_i(\lambda_j)-S_i(\lambda_k)$ is described in the following

\begin{Lem} \label{lem:raction}
Let the action of the group $ \mathbb{R}^+ \times \mathbb{Z}_5$ be defined as above.
Then
\begin{itemize}
\item[(i)] $(x,m)$ leaves the graph \cal{S} invariant, but for a shift of the labels $k \to k+m$
of the external vertices.
\item[(ii)]$(x,m)[S_i(\lambda_j)- S_i(\lambda_k)]= x^{\frac{5}{2}}
\left( S_i(\lambda_j)- S_i(\lambda_k) \right) \;.$ 
\item[(iii)]$(x,m)[\rho_l^k]= x^{-\frac{5}{2}} \rho_l^k \; .$
\end{itemize}

\begin{proof}
The proof of (i) and (ii) follows from the following equality
$$\sqrt{V(\lambda; \Omega^{2k} x^2 a, \Omega^{3k} x^3 b)}d\lambda = x^{\frac{5}{2}}\sqrt{V(\lambda';a,b)} d \lambda' \;, \quad  \lambda=x \lambda' \,.$$

The proof of point (iii) follows from a similar scaling law of the
1-form $\alpha(\lambda) d\!\lambda$ (see equation (\ref{def:alpha}) in Appendix \ref{app:A}).
\end{proof}

\end{Lem}

Due to Lemma \ref{lem:raction}(iii),  $\e=\av{\frac{a}{b}} $ plays
the role of the small parameter. Indeed, along any orbit of the
action of the group $\mathbb{R} \times \mathbb{Z}_5$, all the (finite) relative errors
go to zero uniformly as  $\av{\frac{a}{b}} \to 0$.

Since all the relevant information is encoded in the quotient of the space of cubic potentials
with respect to the group action, we define the following change of variable

\begin{eqnarray}\label{def:mu,nu}
 \nu(a,b) = \frac{b}{a} \;, \qquad
\mu(a,b) = \frac{b^2}{a^3} \quad.
\end{eqnarray}
The induced action on these coordinates is simple, namely
$$(x,m)[\mu(a,b)]=\mu(a,b) \mbox{ and } (x,m)[\nu(a,b)]= \Omega^m x \, \nu(a,b) \; .$$
Moreover, the orbit of the set $\left\lbrace (\nu,\mu) \in \mathbb{C}^2 \mbox{ s.t. }
\av{\,\nu\,}=1, \av{\,\arg{\nu}\,} <\frac{\pi}{5} \,,\;
\mu\neq 0 \right\rbrace$ is a dense open subset of the space of cubic potentials.

\section{Poles of Integr\'ale Tritronqu\'ee}

From Lemma \ref{lem:exclusion} and the results of the computations in Section \ref{par:values},
it follows that equation (\ref{eq:schr}) admits in WKB approximation the simultaneous solutions
of the two quantization conditions $w_{\pm 1}=w_{\mp2}$ only if the Stokes complex is of type $(320)$.
In particular, after our calculations the poles of the int\'egrale tritronqu\'ee are related,
in WKB approximation, to the solutions of the
\textit{Bohr-Sommerfeld-Boutroux} system
(\ref{eq:bohr1},\ref{eq:bohr2}).

\begin{figure}[htbp]
\begin{center}
\input{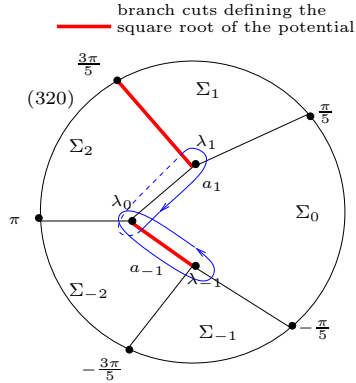}

\end{center}

\caption{Riemann surface $\mu^2=V(\lambda;a,b)$}
\label{figure:bohr}
\end{figure}

We rewrite this system in the following equivalent form:

\begin{eqnarray}\nonumber
 \oint_{a_1}\sqrt{V(\lambda;a,b)}d\lambda &=& i \pi (n-\frac{1}{2}) \\ \label{eq:boutroux} \\ \nonumber
\oint_{a_{-1}}\sqrt{V(\lambda;a,b)}d\lambda &=& -i \pi (m-\frac{1}{2})
\end{eqnarray}

where $m,n$ are positive natural numbers and the paths of integration are shown in figure \ref{figure:bohr}.

System (\ref{eq:boutroux}) is studied in detail in  \cite{kk} where the following lemma
is proven.

\begin{Lem}
 If a polynomial $V(\lambda)= 4 \lambda^3 - 2 a \lambda -28 b$ satisfies
the system (\ref{eq:boutroux}) then
$|\arg{ a}| > \frac{4 \pi}{5}$.
\end{Lem}
The Lemma above should be compared to the conjecture, to which we referred in the introduction.

\paragraph{Real Poles}

We compute all the real solutions of system (\ref{eq:boutroux})
and compare them with some numerical results from
\cite{joshi}. We note that the accuracy of the WKB method is astonishing also
for small $a$ and $b$ (see
Table 2 below).

For the subset of real potentials, we have $$ \oint_{a_1}\sqrt{V(\lambda;a,b)}d\lambda= \overline{\oint_{a_{-1}}\sqrt{V(\lambda;a,b)}d\lambda} \quad,$$
where $\overline{~^{~}}$ stands for complex conjugation.

Therefore system (\ref{eq:boutroux}) reduces to one equation
and the real poles of tritronqu\'ee are characterized, in WKB approximation,
by one natural number.

\begin{Lem}
Let $\mu$ and $\nu$ be defined as in equation (\ref{def:mu,nu}). Then
the real polynomials whose Stokes complex is of type $(320)$ are the orbit of a single point of
the $\mathbb{R}^+$ action, characterized by $\mu^* \cong - 3158,92$ and $\nu>0$.

Moreover, if a real polynomial $V(\lambda; a,b)$ satisfies the Bohr-Sommerfeld-Boutroux conditions
(\ref{eq:boutroux}), then
\begin{equation}\label{eq:realpoles}
a= a^*(n-\frac{1}{2})^{\frac{4}{5}} , \; b= b^* (n-\frac{1}{2})^{\frac{6}{5}}
\end{equation}
 for some $n \in \mathbb{N}^*$ and $a^* \cong -4,0874 , \; b^* \cong -0,1470$.

\begin{proof}
All real cubic potentials whose Stokes complex has type $(320)$ have one real turning point $\lambda_0$
and two complex conjugate turning points $\lambda_{\pm 1}$. For the subset
of real potentials with two complex conjugate turning points the cycles $a_{\pm1}$,
as shown in Figure \ref{figure:bohr}, are unambiguously defined.
From the Classification Theorem,
it follows that a Stokes complex has type $(320)$ if and only if
${\rm Re}\oint_{a_1}\sqrt{V(\lambda;a,b)}d\lambda=0,
{\rm Im}\oint_{a_1}\sqrt{V(\lambda;a,b)}d\lambda \neq 0 $.
Since these conditions are invariant under
the $\mathbb{R}^+$ action, if they are satisfied for a point of the space of cubic potentials,
then they are satisfied on all its orbit. Moreover, it is easily seen that this orbit exists and is unique.
With the help of a software of numeric calculus, we characterized numerically the orbit. Afterthat,
using the scaling law in Lemma \ref{lem:raction}(ii) we calculated all the real solutions of
the Bohr-Sommerfeld-Boutroux system.
\end{proof}

\end{Lem}
To the best of our knowledge, the asymptotic for the $b$ coefficients has never been given.

In the paper \cite{joshi}, the authors showed that the int\'egrale tritronqu\'ee
has no poles on the real positive axis.
The real poles are a decreasing sequence of negative numbers $a_n$ and
some of them are evaluated numerically
in the same paper.

In Table 2, we compare the first two real
solutions to system (\ref{eq:boutroux}) with the numerical evaluation of
the first two poles of the int\'egrale tritronqu\'ee.

\begin{table}
\begin{tabular}{|l|l|l|l|} \hline

& WKB & Numeric & Error \% \\ \hline
$a_1$ & $ -2,34 $ & $ -2,38  $ & $ 1,5 $ \\ \hline
$b_1$ & $-0,064 $ & $-0,062 $ & $2$ \\ \hline
$\mu_1$ & $-3158$ & $-3510$ & $10$ \\ \hline
$a_2$ & $-5,65$ & $-5,66$ & $0,2$ \\ \hline
$b_2$ & $-0,23$ & unknown & unknown \\ \hline
$\mu_2$ & $-3158$ & unknown & unknown \\ \hline
\end{tabular}
\label{table:comparison}
\caption{Comparison between numerical and WKB evaluation of the first two real poles
of the int\'egrale tritronqu\'ee.}
\end{table}

\section{Concluding Remarks}
We have studied the distribution of the poles of solutions to the Painlev\'e first equation
using the theory of the cubic anharmonic oscillator. We have applied a suitable
version of the complex WKB method
to analyze the distribution of poles of the int\'egrale tritronqu\'ee.

In subsequent publications we plan to pursue our study of poles of P-I transcendents in different
directions.

In particular, we want to use the Nevanllina theory of the branched coverings of the sphere
to complete the analysis of the poles of the int\'egrale tritronqu\'ee, showing
that the developed WKB method yields a complete qualitative picture and efficient quantitative
estimates of the distribution of the poles.

Since the quantization condition $\sigma_0=0$ characterizes the monodromy data of a family of
special solutions of P-I, called \textit{int\'egrale tronqu\'ee} (see \cite{kapaev}),
these solutions are strictly related to the spectral theory of PT symmetric
anharmonic oscillators and to functional
equations of Bethe Ansatz type (for what concerns the PT symmetric
anharmonic oscillators and the Bethe Ansatz equations,
see  \cite{tateo} and references therein).
We are going to investigate the consequences of this relation
in a subsequent publication.

\begin{Rem*}
After the main computations of the present paper had been completed, the author learned from B. Dubrovin
about the results of V. Novokshenov presented at the conference NEEDS09 (May 2009).
Novokshenov studied WKB solutions to the Schr\"odinger equation (\ref{eq:schr}) with $b=0$
and their connections to the distributions of poles of certain particular solutions to the
Painlev\'e-I equation, including int\'egrales tronqu\'ee and int\'egrale tritronqu\'ee.
\end{Rem*}

\begin{Rem*}
In the paper \cite{piwkb2} written after this paper was published, the author proves
that eventually (for big enough $n$ and $m$) around any solution of the Bohr-Sommerfeld-Boutroux system
(\ref{eq:boutroux}) there is one and only one pole of the int\'egrale tritronqu\'ee. Moreover the distance
between a pole and its approximation vanishes asymptotically.
\end{Rem*}

\appendix

\section{Appendix} \label{app:A}

The aim of this appendix is to prove  Theorem \ref{thm:birkh}. Our approach is similar to the 
approach of Fedoryuk \cite{fedoryuk}.

Notations are as in sections \ref{section:complexes} and  \ref{section:wkb}, except for 
$\infty_k$.
In what follows, we suppose to have fixed a certain cubic potential $V(\lambda;a,b)$
and a maximal domain $D_k$. To simplify the notation we write $V(\lambda)$ instead of
$V(\lambda;a,b)$.

\subsection{Gauge Transform to an L-Diagonal System}

The strategy is to find a suitable gauge transform of equation (\ref{eq:schr}) such that
for large $\lambda$ it simplifies.
We rewrite the Schr\"odinger equation 
 
\begin{eqnarray}\label{eqn:scalar}
-\psi''(\lambda)+V(\lambda) \psi(\lambda)=0 \, \; ,
\end{eqnarray}

in first order form:
\begin{eqnarray}
\Psi'(\lambda)=E(\lambda)\Psi(\lambda) \nonumber \; , \\ \label{eqn:matricial}
E(\lambda)=\left(  \begin{matrix}
0 & 1 \\ 
V(\lambda) & 0
\end{matrix}\right) \; .
\end{eqnarray}

\begin{Lem}[Fedoryuk]\label{lem:gauge}

In $D_k$

\begin{itemize}
\item[(i)] the gauge transform 
\begin{eqnarray}
Y(\lambda) &=&A(\lambda)U(\lambda) \nonumber \; , \\ \label{eqn:gaugetransform}
A(\lambda) &=& j_k(\lambda) \left(  \begin{matrix}
1 & 1 \\ 
\sqrt{V(\lambda)}-\frac{V'(\lambda)}{4V(\lambda)} & -\sqrt{V(\lambda)}-\frac{V'(\lambda)}{4V(\lambda)} 
\end{matrix}\right) \; ,
\end{eqnarray}
is non singular and
\item[(ii)] the system (\ref{eqn:matricial}) is transformed into the following one

\begin{eqnarray}
U'(\lambda) &=&F(\lambda)U(\lambda) = \left( A^{-1}EA - A^{-1}A' \right) U \; ,\nonumber\\
 \label{eqn:transformed}
F(\lambda) &=& 2 \sqrt{V(\lambda)}\left(  \begin{matrix}
1 & 0 \\ 
0 & 0 
\end{matrix}\right) + \alpha(z) \left(  \begin{matrix}
1 & 1 \\ 
-1 & -1 
\end{matrix}\right) \; , \\
\label{def:alpha}
\alpha(\lambda) &=& \frac{1}{32 \sqrt{V(\lambda)}^5}(4 V(\lambda) V''(\lambda) - 5 V'^2(\lambda)) \; .
\end{eqnarray}

\end{itemize}
\begin{proof}

\begin{itemize}
\item[(i)] Indeed $\det A(\lambda)= 2 j_k^2(\lambda)\sqrt{V(\lambda)} \neq 0 , \; \forall \lambda \in D_k$, 
by construction of $j_k$ and $D_k$.
\item[(ii)] It is proven by a simple calculation.
\end{itemize}
\end{proof}
\end{Lem}

\subsection{Some Technical Lemmas}
Before we can begin the proof of Theorem \ref{thm:birkh}, we have to introduce the compactification
of $D_k$ and the preparatory Lemmas \ref{lem:prep1} and \ref{lem:prep2}.

\paragraph{Compactification of $D_K$}

Since $D_k$ is simply connected, it is conformally equivalent to the interior of the unit disk $D$.
We denote $U$  the uniformisation map, $U: D \to D_k$.

By construction, the boundary of $D_k$ is the union of $n$ free Jordan curves,
all intersecting at $\infty$. Here $n$ is
equal to the number of 
sectors $\Sigma_l$ such that $ \Sigma_l\leftrightarrows \Sigma_k$ minus 2.

Due to an extension of Carath\'eodory's Theorem (\cite{carat}, \S 134-138),
the map $U$ extends to a continuous map from the closure of the unit circle
to the closure of $D_k$. The map is injective on the closure of $D$ minus 
the $n$ counterimages of $\infty$.
Hence, the uniformisation map realizes a $n$ point compactification of $D_k$,
that we call $\overline{D_k}$.
In $\overline{D_k}$ there are $n$ point at $\infty$.
We denote $\infty_k$ the point at $\infty$ belonging to the closure of
$U\left( \Sigma_{k-1} \cup \Sigma_{k} \cup \Sigma_{k+1} \right)$. Moreover, for
$\lambda = k +2 \mbox{ or } \lambda=k-2$, if $\Sigma_{l} \leftrightarrows \Sigma_k$
we denote $\infty_{l}$ the point at $\infty$ belonging to the closure of $U(\Sigma_{l})$.

\begin{Def}\label{def:K}
Let $H$ be the space of  function holomorphic in  $D_k$ and continuous in $\overline{D_k}$.
$H$ endowed with the sup norm is a Banach space ($H$, $\lVert \cdot \rVert_H$).

Let $\Gamma(\lambda), \lambda \in \overline{D_k}-\infty_k$ be the set of injective piecewise
differentiable curves $\gamma:[0,1] \to \overline{D_k}$,
such that
\begin{enumerate}
 \item $\gamma(0)=\lambda$, $\gamma(1)=\infty_k$,
 \item $ReS_k(\gamma(0),\gamma(t))$ is eventually non decreasing,
 \item there is an $\e >0$ such that eventually
 $\left| \arg{\gamma(t)} -\frac{2\pi k}{5}\right| < \frac{\pi}{5} -\varepsilon$,
 \item the length of the curve restricted to $[0,T]$
  is $O \left(\left|\gamma(T)\right| \right), \mbox{as } t \to 1$.
\end{enumerate}

Let $\tilde{\Gamma}(\lambda)$ be the subset
of $\Gamma(\lambda)$ of the paths along which $ReS_k(\gamma(0), \gamma(t))$
is non decreasing.

Let $K_1: H \to H$ and $K_2: H \to H$ be defined (for the moment formally) 

\begin{eqnarray}\label{def:K1}
K_1 [h](\lambda) &=& -\int_{\gamma \in \Gamma(\lambda)} e^{2S_k(\mu,\lambda)}
\alpha(\mu) h(\mu) d\mu \;,\\ \label{def:K2}
K_2 [h](\lambda)  &=& \int_{\gamma \in \Gamma(\lambda)} \alpha(\mu) h(\mu) d\mu \; .
\end{eqnarray}

Let $\rho: \overline{D_k} \to \overline{D_k}$:
 \begin{equation*}
\rho(\lambda)= \left\lbrace 
             \begin{aligned}
               \inf_{\gamma\in \tilde{\Gamma}(\lambda)}
\int_{0}^1 \left| \alpha(\gamma(t))\frac{d\gamma(t)}{dt} \right| dt, \; \;
\mbox{if} \;\; \lambda \neq \infty_k \\
 0 \qquad \qquad , \;\; \mbox{if} \;\; \lambda=\infty_k \; .
             \end{aligned}
\right.  
\end{equation*}
\end{Def}

\begin{Rem*}
Since along rays of fixed argument $\varphi$, with  $\left| \varphi -\frac{2\pi k}{5}\right| < \frac{\pi}{5} -\varepsilon$, $ReS_k$ is eventually increasing, there are paths satisfying point (1) through (4) of the above definition. Moreover, by construction of $D_k$, $\tilde{\Gamma}(\lambda)$ is non empty for any $\lambda$.
\end{Rem*}

Before beginning the proof of the theorem, we need two preparatory lemmas.

\begin{Lem}\label{lem:prep1}
Fix $\varepsilon >0$, an angle $\left| \arg{\varphi} -\frac{2\pi l}{5}\right| < \frac{\pi}{5} -\varepsilon $,
and let $\Omega =\Sigma_l \cap \left\lbrace \lambda \in \mathbb{C},
\av{\lambda - \frac{2 \pi l}{5} }< \frac{\pi}{5} - \e \right\rbrace$.
Denote $i(R) = i_{R e^{i \varphi}} \cap \Omega , R \in \mathbb{R}^+$,
and let $L(R)$ be the length with respect to the euclidean metric of $i(R)$.
Then $L(R)=O(R)$ and
$\inf_{\lambda' \in i(R)}\left|\lambda'\right|=O(R)$.

Let $r$ be any level curve of $S_l(\lambda^*,\cdot)$ asymptotic to the ray of argument $\frac{2 \pi l}{5}$, $\Omega(R)=
\left\lbrace \lambda \in \Omega, ReS_l(\lambda, R e^{i \varphi}) \geq 0 \right\rbrace$,
and $M(R)$ be the length
of $r \cap \Omega(R)$. Then $M(R)=O(R)$.

\begin{proof}
\cite{strebel}, chapter 3. 
\end{proof}

\end{Lem}

\begin{Lem}\label{lem:prep2}
\begin{itemize}
 \item [(i)] $\rho$ is a continuous function.
\item[(ii)] $K_1$ and $K_2$ are well-defined bounded operator. In particular
\begin{equation}\label{eq:estimate}
 \left| K_i [h] (\lambda) \right| \leq  \rho(\lambda) \lVert h \rVert_{H} \, , \; i=1,2
\end{equation}
\item[(iii)]$K_2[h](\infty_k)=K_1[h](\infty_k)=K_1(\infty_{k \pm 2})=0, \forall h \in H$

\end{itemize}

\begin{proof}
(i)Since $\alpha(\lambda)d\lambda= O(\left| \lambda \right|^{-\frac{7}{2}})$, then
$\alpha(\lambda)d\lambda$ is integrable along any curve $\gamma \in \tilde{\Gamma}(\lambda)$.
Therefore $\rho$ is a continuous function on $\overline{D_k}$.

(ii)We first prove that (a) $K_i[h](\lambda)$  does not depend on the integration path 
for any $\lambda \in \overline{D_k}$ minus the points at infinity. A result that easily implies that
$K_i[h](\cdot)$ is an analytic function on $D_k$, continuous on $\lambda \in \overline{D_k}$
minus the points at $\infty$.
We then prove (b) the estimates (\ref{eq:estimate}) and (c) the existence of the limits
$K_i[h](\infty_l), \; l = \infty_k, \infty_{k \pm 2}$.

To simplify the notation, we prove the theorem for the operator $K_1$. The 
proof for $K_2$ is almost identical .

(a)Let $\gamma_a, \gamma_b \in \Gamma(\lambda)$. The curve $i_{\gamma_a(T)}$,
where $T= 1 - \e$ for some small $\e >0$ intersect $\gamma_b$
at some $\gamma_b(T')$. Therefore we can decompose $-\gamma_b \circ \gamma_a$
into two different paths with the help of a segment of $i_{\gamma_a(T)}$,
$\int_{-\gamma_b \circ \gamma_a}e^{2S_k(\mu,\lambda)} \alpha(\mu) h(\mu) d\mu
=\int_{\gamma_1}+\int_{\gamma_2}e^{2S_k(\mu,\lambda)} \alpha(\mu) h(\mu) d\mu$.
One path $\gamma_1$ is the loop based at $\lambda$ and the other $\gamma_2$
is the loop based at $\infty_k$.
Since $\gamma_1 \subset D_k$, then $\int_{\gamma_1}e^{2S_k(\mu,\lambda)} \alpha(\mu) h(\mu) d\mu=0$.
Along $\gamma_2$, $e^{2S_k(\gamma_2(t),\lambda)} \leq 1$ therefore 
the integrand can be estimated just by $\left| \alpha(\gamma_2(t)) \right|$. Due to lemma
\ref{lem:prep1}, $\int_{\gamma_2}\av{\alpha(\mu) h(\mu) d\mu}=O(\av{\gamma_a(T)}^{-\frac{5}{2}})$.
Since $\e$ is arbitrary, then $K_1[h](\lambda)$ does not depend on the integration path.

(b)Clearly for any path $\gamma \in \tilde{\Gamma}(\lambda)$,
$\left| K_1 [h] (\lambda) \right| \leq \int_{0}^1 \av{\alpha(\lambda) h(\lambda)
d\!\lambda} dt$. Since $K_1[h](\lambda)$ does not
depend on $\gamma$, then  estimate (\ref{eq:estimate}) follows.

(c) Let $\lambda_n$ be a sequence converging to $\infty_{l}, l=k+2 \mbox{ or } l=k-2$;
without losing any generality we suppose that
the sequence is ordered such that $ReS_k(\lambda_n) \leq ReS_k(\lambda_{n+1})$.
Fix a curve $r$, as defined in Lemma \ref{lem:prep1}.
By construction of $D_k$, it is always possible to connect two points $\lambda_{n}$ and $\lambda_{n+m}$
with a union of segments of the curves $i_{\lambda_n}, i_{\lambda_{n+m}}$ and of $r$.
We denote by $\gamma$ the union of this three segment.
By construction of $D_k$ (see Subsection \ref{sub:domains} (iii)), there exists $\e\!>\!0$
such that $\av{\arg{\lambda_n} - \frac{2 \pi l}{5}} < \frac{\pi}{5} - \e , \forall n$. Therefore,
due to Lemma \ref{lem:prep1},  $\gamma$ has length
of order $\left| \lambda_{n} \right| + \left| \lambda_{n+m} \right|$.
Hence $\left| K_1[h](\lambda_n)-K_1[h](\lambda_{n+m}) \right| \leq
\int_{\gamma}\left|h(\lambda) \alpha(\lambda)  d\lambda \right| = O(\av{\lambda_n}^{-\frac{5}{2}})$.
Then $K_1[h](\lambda_n)$ is a Cauchy sequence
and the limit is well defined.

We now prove that this limit is zero by calculating
it along a fixed ray $\lambda=x e^{i \varphi}$ inside $\Sigma_{k \pm 2}$. Let us fix
a point $x^*$ on this ray in such a way that the function $ReS_k(x^*,x)$ is monotone decreasing
in the interval $[x^*,+\infty[$. Along the ray we have
$$K_1[h](x)= - \frac{\int_x^{x^*} e^{2S_k(y,x^*)}\alpha(y) h(y) dy + g(x^*)  }{ e^{2S_k(x^*,x)}}\;, $$
where $g(x^*)$ is a constant, namely $\int_{\gamma \in \Gamma(x^*)} e^{2S_k(\mu,x^*)} \alpha(\mu) h(\mu) d\mu$.
Hence $\lim_{x \to \infty}K_1[h](x)=\lim_{x \to \infty}\frac{\alpha(x)h(x)}{\sqrt{V(x)}}=0$.

With similar methods the reader can prove that the limit  $K_1[h](\infty_k)$ exists
and is zero.

\end{proof}

\end{Lem}

We are now ready to prove Theorem \ref{thm:birkh}.

\begin{Thm}
Extend the WKB function $j_k$ to $D_k$.
There exists a unique solution $\psi_k$ of (\ref{eq:schr}) such that for all $\lambda \in D_k$
\begin{eqnarray*}
\left| \frac{\psi_k(\lambda)}{j_k(\lambda)} -1 \right| &\leq& g(\lambda)(e^{2\rho(\lambda)}-1) \;, \\
\left| \frac{\psi_k'(\lambda)}{j_k(\lambda)\sqrt{V(\lambda)}} +1 \right|
&\leq&  \left|\frac{V'(\lambda)}{4V(\lambda)^{\frac{3}{2}}}\right|+
(1+\left|\frac{V'(\lambda)}{4V(\lambda)^{\frac{3}{2}}}\right|)g(\lambda)( e^{2 \rho(\lambda)} -1) \; , 
\end{eqnarray*}

where $g(\lambda)$ is a positive function, $g(\lambda) \leq 1$ and $g(\infty_{k\pm 2})=\frac{1}{2}$.

\begin{proof}

We seek a particular solution to the linear system (\ref{eqn:transformed}) via 
successive approximation.

If $U(\lambda) = U^{(1)} \oplus U^{(2)} \in H \oplus H$
satisfies the following integral equation of Volterra type

\begin{eqnarray}\nonumber
U(\lambda) &=& U_0 + K[U](\lambda)  \; , \; U_0 \equiv \left( \begin{matrix}
0 \\ 
1
\end{matrix} \right) \; ,  \\ \label{eqn:K}
K[U](\lambda) &=&  \left( \begin{matrix}
K_1[U^{(1)}+U^{(2)}](\lambda)\\ 
K_2[U^{(1)}+U^{(2)}](\lambda)
\end{matrix} \right) \; ,
\end{eqnarray}

then $U(\lambda)$ restricted to $D_k$ satisfies (\ref{eqn:transformed}).

We define the the \textit{Neumann series} as follows

\begin{equation}\label{def:neumann}
U_{n+1}=U^0+K[U_n]  \, , U_{n+1}= \sum_{i=0}^{n+1}K^i[U^0] \; .
\end{equation}

More explicitly,

\begin{eqnarray*}
K^n[U_0](\lambda)= \left( \int_{\lambda}^{\infty_k} d\mu_1 \int_{\mu_1}^{\infty_k} d\mu_2 \dots \int_{\mu_{n-1}}^{\infty_k} d\mu_n  \begin{matrix}
 -e^{2S(\mu_1,z)} \alpha(\mu_1) &\times &\\ 
\alpha(\mu_1) &\times &
\end{matrix} \right. \\
\left.
\begin{matrix}
 \alpha(\mu_2)(1-e^{2S(\mu_2,\mu_1)}) \dots \alpha(\mu_n)(1-e^{2S(\mu_{n},\mu_{n-1})}) \\ 
 \alpha(\mu_2)(1-e^{2S(\mu_2,\mu_1)}) \dots \alpha(\mu_n)(1-e^{2S(\mu_{n},\mu_{n-1})})
\end{matrix} \right) \; .
\end{eqnarray*}

Here the integration path $\gamma$ belong to $\Gamma(\lambda)$.
For any $\gamma \in \tilde{\Gamma}(\lambda)$ and any $n \geq 1$

\begin{eqnarray*}
\left|\! K^n[U_0]^{(i)}(\lambda) \!\right| \! \leq\!  \frac{1}{2} \int_{\lambda}^{\infty_k} \! \int_{\mu_1}^{\infty_k}\!
\dots \! \int_{\mu_{n-1}}^{\infty_k} \prod_{i=1}^{n} |2\alpha(\mu_i)d\!\mu_i|\! =\!
\frac{2^{n-1}}{n!} \! \left(\! \int_{\gamma} d\mu_1 |\alpha(\mu_1)| \!\right)^n ,
\end{eqnarray*}

where $K^n[U_0]^{(i)}$ is the i-th component of $K^n[U_0]$. Hence

\begin{equation}\label{ineq:Kn}
| K^n[U]_i(\lambda) | \leq \frac{1}{2} \frac{1}{n!} \left(2 \rho(\lambda) \right)^n
\end{equation}

Thus the sequence $U^n$ converges in $H$ and is a solution to (\ref{eqn:K}); call $U$ its limit.
Due to Lemma \ref{lem:prep2}, $U^{(1)}(\infty_{k \pm 2})=0$.

Let $\Psi_k$ be the solution to (\ref{eqn:matricial}) whose gauge transform is $U$ restricted to $D_k$;
The first component $\psi_k$ of $\Psi_k$ satisfies equation (\ref{eqn:scalar}).

From the gauge transform (\ref{eqn:gaugetransform}), we obtain 

\begin{eqnarray*}
\frac{\psi_k(\lambda)}{j_k(\lambda)} -1 &=&U_1(\lambda)+U_2(\lambda)-1 \; ,\\
\frac{\psi_k'(\lambda)}{j_k(\lambda)\sqrt{V(\lambda)}} +1 &=&U_1(\lambda)(1-\frac{V'(\lambda)}{4V(\lambda)^{\frac{3}{2}}})-
(U_2(\lambda)-1)(1+\frac{V'(\lambda)}{4V(\lambda)^{\frac{3}{2}}})+\\
&-&\frac{V'(\lambda)}{4V(\lambda)^{\frac{3}{2}}} \; ,
\end{eqnarray*}

The thesis follows from these formulas, inequality (\ref{ineq:Kn}) and from the fact that
$U_1(\infty_{k \pm 2})=0$.

\end{proof}
\end{Thm}

\begin{Rem*}
The solution $\psi_k(\lambda)$ of equation (\ref{eq:schr}) described in Theorem \ref{thm:birkh}
may be extended from $D_k$ to the whole complex plane,
since the equation is linear with entire coefficients. The continuation
is constructed in the following Corollary.
\end{Rem*}

\begin{Cor}\label{cor:continuation}
For any $\lambda \in \mathbb{C}$, $\lambda$ not a turning point,
we define $\Gamma(\lambda)$ as in Definition \ref{def:K}.
Fixed any $\gamma \in \Gamma(\lambda)$ and $h$ a continuous function on $\gamma$,
we define the functionals $K_i[h](\lambda)$ as in equations 
(\ref{def:K1}) and (\ref{def:K2}). We define the Neumann series as in equations (\ref{eqn:K}) and (\ref{def:neumann}), and we continue $j_k$ along $\gamma$. 

Then then Neumann series converges and we call $U^{(1)}(\lambda)$ and $U^{(2)}(\lambda)$
the first and second component of its limit.

Moreover, $\psi_k(\lambda)= \left( U^{(1)}(\lambda)+
U^{(2)}(\lambda) \right) j_k(\lambda)$
solves equation (\ref{eq:schr}) and for any $\e >0$
$$ 
\lim_{\av{\lambda} \to \infty \, ,\; \av{arg{\lambda}-\frac{2 \pi k}{5}}<\frac{3 \pi}{5} - \e}\left( U^{(1)}(\lambda)+U^{(2)}(\lambda) \right)=1
$$

\end{Cor}

The reader should notice that if $\lambda \notin D_k$, then $\tilde{\Gamma}(\lambda)$ is empty and
we cannot estimate $\frac{\psi_k(\lambda)}{j_k(\lambda)}$.

\section{Appendix}\label{app:B}

The aim of this Appendix is to prove Theorem \ref{them:toCubic} (ii) and (iii).
The notation is, if not otherwise stated, as in the previous sections of the paper.

Next to a pole $z=a$ of a solution $y(z)$ of P-I, equation (\ref{sys:1}) becomes
meaningless. To get rid of this singularity we perform a gauge transform of (\ref{sys:1})
such that the gauge-transformed equation has full meaning in the limit.
In what follows, we suppose that $z$ belongs to a punctured neighborhood of
$a$, where $y(z)$ is holomorphic.

\subsection*{A gauge transform}
Let $z$ be a fixed regular value of $y(z)$. Let $\overrightarrow{\Phi}(\lambda,z) =
G(\lambda, z) \overrightarrow{\Psi}(\lambda,z)$,

\begin{eqnarray}\label{toscalar}
G(\lambda,z) = \left(
\begin{matrix}
\frac{y'(z)+\frac{1}{2(\lambda-y(z))}}{\sqrt{2(\lambda-y(z))}} & \frac{1}{\sqrt{2(\lambda-y(z))}} \\
 & \\
\sqrt{2(\lambda-y(z))} & 0
\end{matrix}
  \right) \; .
\end{eqnarray}

Then  $\overrightarrow{\Phi}(\lambda;z)$ satisfies (\ref{sys:1}) if and only if
$\overrightarrow{\Psi}(\lambda;z)$ satisfies 
the following equation

\begin{eqnarray*}
\Psi_{\lambda}(\lambda,z)= \left(
\begin{matrix}
0 & 1 \\ 
Q(\lambda;z) & 0
\end{matrix}
  \right)
 \Psi(\lambda,z)
\end{eqnarray*}

where 

\begin{equation}\label{perturbedpot}
 Q(\lambda;z)= 4 \lambda^3 - 2 \lambda z + 2 z y(z) - 4 y^3(z) + y'^2(z)+  \frac{y'(z)}{\lambda - y(z)} 
+\frac{3}{4(\lambda -y(z))^2}
\end{equation}

We denote $\psi$ the first component of $\overrightarrow{\Psi}$.
The equation for $\overrightarrow{\Psi}$ is equivalent to the following second order scalar
equation for $\psi$

\begin{equation}\label{eq:perturbedscalar}
 \psi_{\lambda \,\lambda}(\lambda,z)= Q(\lambda;z) \psi(\lambda,z) 
\end{equation}

We summarize some property of the perturbed potential, which can be easily verified using the expansion
(\ref{laurent}).

\begin{Lem}\label{lem:pertpotential}
Let $\varepsilon^2=\frac{1}{y(z)} = (z-a)^2 + O((z-a)^6)$  then
\begin{itemize}
\item[(i)] $Q(\lambda; z)$ has a double pole at $\lambda=\frac{1}{\varepsilon^2}$. It is
an apparent fuchsian singularity for equation (\ref{eq:perturbedscalar}): the local monodromy around
it is $-1$.
\item[(ii)] $Q(\lambda; z)$ has two simple zeros at $\lambda=\frac{1}{\varepsilon^2}+ O(\varepsilon^2)$
\item[(iii)] $Q(\lambda; z)= 4 \lambda^3 -2 (a+\e) \lambda - 28 b + O(\varepsilon)-
\frac{ 2 \lambda \varepsilon^{-1}}{\lambda - \varepsilon^{-2}} + \frac{3}{4 (\lambda - \varepsilon^{-2})^2}$,
where $O(\varepsilon)$ does not depend on $\lambda$.
\end{itemize}

\end{Lem}

Equation (\ref{eq:perturbedscalar}) is a perturbation of the cubic Schr\"odinger equation (\ref{eq:schr})
and the asymptotic behaviours of solutions to the two equations are very similar.
Indeed the local picture around the point at $\infty$ depends only on the terms $4 \lambda^3$
and $-2z\lambda^3$.

More precisely, the equivalent of Corollary \ref{cor:continuation} in Appendix \ref{app:A} is valid also
for the perturbed Schr\"odinger equation.

\begin{Def}\label{def:perturbedwkb}
For any $z$, define a cut from $\lambda=\frac{1}{(z-a)^2}$ to $\infty$ such that
it eventually does not belong to the the angular sector
$\left| \arg{\lambda} -\frac{2 \pi k}{5} \right| \leq \frac{3\pi}{5}$.

Fix $\lambda^*$ in the cut plane. $S_k(\lambda;z)=\int_{\lambda^*}^{\lambda}\sqrt{Q(\mu;z)}d\mu$
is well-defined for $\av{\arg{\lambda}-\frac{2 k \pi}{5}} < \frac{3\pi}{5}$ and $\lambda>>0$.
Here the branch of $\sqrt{Q}$ is chosen such that $ReS_k(\lambda)\to + \infty$ as $\av{\lambda} \to \infty \,,\;\av{\arg{\lambda}-\frac{2 \pi k}{5}}< \frac{\pi}{5}-\e$. We define $j_k(\lambda;z)$ as in equation (\ref{def:jk}) and $\alpha(\lambda;z)$ as in equation (\ref{def:alpha}), but replacing
$V(\lambda)$ with $Q(\lambda;z)$.

For any $\lambda$ in the cut plane, let $\Gamma(\lambda)$ be the set of piecewise differentiable curves
$\gamma:[0,1]$ to the cut plane, $\gamma(0)=\lambda$, $\gamma(1)=\infty$, satisfying
properties (2)(3) and (4) of Definition \ref{def:K}.

For any $\gamma \in \Gamma(\lambda)$, let $H$ be the Banach space of continuous functions
on $\gamma$ that have a finite limit as $t \to 1$.
Formulae (\ref{def:K1}) and (\ref{def:K2}) define two bounded functionals on $H$.
We call such functionals $K_1(\lambda;z)$ and $K_2(\lambda;z)$.
 \end{Def}

Following the proof of Theorem \ref{thm:birkh},
the reader can prove the following

\begin{Lem}\label{lem:perturbedwkb}
Let $\lambda$ belong to the cut plane, $\lambda$ not a zero of $Q(\cdot;z)$.
Fixed any $\gamma \in \Gamma(\lambda)$, we define the Neumann series as in equations (\ref{eqn:K}) and (\ref{def:neumann}), and we continue $j_k$ along $\gamma$.

Then the Neumann series converges and
$\psi_k(\lambda)= \left( U_1(\lambda)+U_2(\lambda) \right) j_k(\lambda)$
solves equation (\ref{eq:perturbedscalar}). Moreover, for any $\e >0$
$$ 
\lim_{\av{\lambda} \to \infty \, ,\; \av{arg{\lambda}-\frac{2 \pi k}{5}}<\frac{3 \pi}{5} - \e}\left( U^{(1)}(\lambda)+U^{(2)}(\lambda) \right)=1
$$
\end{Lem}

\begin{Def}
Let $\tilde{\psi_k}(\lambda,z)$  be the unique solution of equation
(\ref{eq:perturbedscalar}) such that

\begin{equation}\label{scalarasymptotic}
\frac{\tilde{\psi_k}(\lambda,z)} {\lambda^{-\frac{3}{4}} e^{-\frac{4}{5} \lambda^{\frac{5}{2}} +
z\lambda^{\frac{1}{2}}}} \to 1 \, ,  \mbox{ as } \av{\, \lambda \,} \to \infty \,,\,
arg{\lambda}=\frac{2 \pi k}{5} \; .
\end{equation}

Here the branch of $\lambda^{\frac{1}{4}}$ is fixed as $\lambda \to \infty ,\,
\arg{\lambda}=\frac{2 \pi k}{5}$, and there it coincides with the branch chosen
in equation (\ref{asymptotic}).
We define $\tilde{\psi_k}(\lambda,a)$ to be the unique solution of equation (\ref{eq:schr})
with asymptotic (\ref{scalarasymptotic}), where $z=a$.

We denote $\psi_k(\lambda;z)$ the unordered pair $\left\lbrace  \tilde{\psi_k}(\lambda,z),
-\tilde{\psi_k}(\lambda,z) \right\rbrace$.
\end{Def}

\begin{Rem*}
We notice that if the cuts are continuous in $z$, then $\tilde{\psi_k}(\lambda,z)=c(z)\psi_k(\lambda)$,
where $\psi_k(\lambda)$ is the solution constructed in Lemma \ref{lem:perturbedwkb}
and $c(z)$ is a bounded holomorphic function.
\end{Rem*}

\begin{Thm}\label{thm:limit}
$\lim_{z \to a}\psi_k(\lambda, z) = \psi_k(\lambda,a)$, $\forall \lambda \in \mathbb{C}$.

\begin{proof}

Let $\lambda$ be any point in the complex plane
which is not a zero of $V(\lambda;a,b)$.
For any sequence $\e_n$ converging to zero,
we choose two fixed rays $r_1$ and $r_2$ of different argument  $\varphi_1$ and $\varphi_2$,
$\av{\varphi_i- \frac{2 k \pi}{5}}<\frac{\pi}{5}$.
We denote $D_{R,\e}$ a disk of radius $R$ with center $\lambda=\frac{1}{e^2}$ and
we split the sequence $\e_n$ into two subsequences $\e^i_n$ such that
$r_i \cap D_{R,\e^i_n} = \emptyset$ for any $n$
big enough.

For $i=1,2$, we choose the cuts defined in Definition \ref{def:perturbedwkb}
in such a way that there exists a differentiable curve $\gamma_i:[0,1] \to \overline{\mathbb{C}}$,
$\gamma_i(0)=\lambda$, $\gamma_i(1)=\infty$ with the following properties: (i)$\gamma_i$ avoids the zeroes
of $Q(\lambda,\e^i_n)$ and a fixed, arbitrarily small, neighborhood of the zeroes of $V(\lambda;a,b)$,
(ii)$\gamma_i$ does not intersect any cut, and
(iv) $\gamma_i$ eventually lies on $r_i$.

The proof of the thesis relies on the following estimates:

\begin{eqnarray}\label{estimates}
&\sup_{\lambda \in \mathbb{C}-D_{R,\e} }&
\av{\lambda^{-\delta}}\av{\,Q(\lambda;a+\e)-V(\lambda;a,b)\,}= O(\e^{2 \delta - 3}) \; ,\\ \nonumber
&\sup_{\lambda \in \mathbb{C}-D_{R,\e} }&
\av{\lambda^{-\delta}}\av{\,Q_{\lambda}(\lambda;a+\e)-
V_{\lambda}(\lambda;a,b)\,}= O(\e^{2 \delta - 3})\; , \\ \nonumber
&\sup_{\lambda \in \mathbb{C}-D_{R,\e} }&
\av{\lambda^{-\delta}}\av{\,Q_{\lambda\,\lambda}(\lambda;a+\e)-
V_{\lambda\,\lambda}(\lambda;a,b)\,}= O(\e^{2 \delta - 3}) \; .
\end{eqnarray}

Due the above estimates it is easily seen that $\gamma_i \in \Gamma(\lambda)$, $\forall \e^i_n$.
Due to Lemma \ref{lem:perturbedwkb} and Corollary \ref{cor:continuation},
to prove the thesis it is sufficient to show that the functionals
$K_1(\lambda; a+\e^i_n)$ and $K_2(\lambda;a+\e^i_n)$
converge in norm to $K_1(\lambda;a)$ and $K_2(\lambda;a)$. Here
$K_i(\lambda;a), \, i=1,2$ are defined as in Corollary \ref{cor:continuation}.
We notice that the norm of the functionals are
just the $L^1(\gamma_i)$ norm of their integral kernels.

We first consider the functionals $K_2(\lambda;a +\e^i_n)$. Due to the above estimates

\begin{equation*}
\lambda^{\frac{7}{2}}\alpha(\mu,\e^i_n) \to \lambda^{\frac{7}{2}}\alpha(\mu), \text{ uniformly on } \gamma_i([0,1])
\text{ as } \; n \to \infty \; .
\end{equation*}

Hence the sequence $\alpha(\mu,\e^i_n)$ converges in norm $L^1(\gamma_i)$
to $\alpha(\mu)$ and  the sequence $K_2(\lambda; e_n^i)$ converges in
operator norm to $K_2(\lambda;a)$.

We consider now the sequence $K_1(\lambda;a + \e^i_n)$.

To prove the convergence of the above sequence of operators, it is sufficient to prove that

\begin{equation*}
e^{S_k(\lambda; a+\e^i_n)-S_k(\mu;a+\e^i_n)} \to e^{S_k(\lambda; 0)-S_k(\mu;0)} \text{ uniformly on } \gamma_i([0,1]) \text{ as } n \to \infty \; .
\end{equation*}

We first note that

\begin{eqnarray*}
& e^{S_k(\lambda; 0)-S_k(\mu;0)}& - e^{S_k(\lambda; a+\e^i_n)-S_k(\mu;a+\e^i_n)}
=  e^{S_k(\lambda; 0)-S_k(\mu;0)} \left( 1- e^{g(\mu; \e)} \right) \; ,\\
& g(\mu, \e) &= {\int_{\lambda, \gamma_i}^{\mu} \frac{Q(\nu, \e)-V(\nu;a,b)}
{\sqrt{Q(\nu,\e)} + \sqrt{P(\nu;a,b)}} \, d \nu} \; .
\end{eqnarray*}

Using estimate (\ref{estimates}), it is easy to show that $g(\mu; \e)= f(\e)O(\mu^{\delta})$,
where $f(\e) \to 0 \text{ as } \e \to 0$ and $0<\delta<<1$. Therefore the difference of
the exponential functions converges uniformly to 0.

\end{proof}
\end{Thm}

We can prove  Theorem \ref{them:toCubic} (ii) and (iii).

Indeed from (\ref{eqn:gaugetransform}),
it is easily seen that (choosing
one of the two branches of the gauge transform)

$$\overrightarrow{\Psi_k}(\lambda;z)=G(\lambda,z)^{-1}\overrightarrow{\Phi_k}(\lambda;z)= 
\frac{1}{\sqrt{2} } \left(
\begin{matrix}
 \tilde{\psi}_k(\lambda,z) \\ 
\tilde{\psi_k'}(\lambda,z)
\end{matrix} \right) \; , $$

if $\av{\, \arg{\lambda}-\frac{2 \pi k}{5}\,} <\frac{3 \pi}{5} \mbox{ and } \av{\lambda}>>0$.

Moreover from (\ref{eqn:gaugetransform}), it follows that 

$$\lim_{z\to a} (z-a) \Phi_k^{(2)}(\lambda,z)= i \sqrt{2}\Psi_k^{(1)}(\lambda;a) \; .$$

Hence Theorem \ref{them:toCubic} (ii) and (iii) follow from Theorem \ref{thm:limit}.

\newpage


\bibliographystyle{alpha}
\newcommand{\etalchar}[1]{$^{#1}$}

\end{document}